\documentclass[1pt,notitlepage,twoside,a4paper]{amsart}

\usepackage{amsmath,amssymb,enumerate}

\usepackage{epsfig,fancyhdr,color}

\usepackage{amssymb}
\usepackage{amsmath,amsthm}    
\usepackage{latexsym}
\usepackage{amscd} 
\usepackage{psfrag}
\usepackage{graphicx} 
\usepackage[latin1]{inputenc}  
\usepackage[all]{xy} 
\usepackage[mathcal]{eucal}

\theoremstyle{definition}

\theoremstyle{remark}

\def\interieur#1{\mathord{\mathop{\kern 0pt #1}\limits^\circ}}

\definecolor{NoteColor}{rgb}{1,0,0}

\title[Physics in Riemann's mathematical papers]{Physics in Riemann's mathematical papers}

\author{Athanase Papadopoulos}
\address{A. Papadopoulos, Institut de Recherche Math{\'e}matique Avanc\'ee,
Universit{\'e} de Strasbourg and CNRS,
7 rue Ren\'e Descartes,
 67084 Strasbourg Cedex, France,  and 
  Brown University, Mathematics Department, 
 151 Thayer Street
Providence, RI 02912, USA.}

 \date{\today}

\makeindex

\begin{document}

\begin{abstract}

Riemann's mathematical papers contain many ideas that arise from physics, and some of them are motivated by problems from physics. In fact, it is not easy to separate Riemann's ideas in mathematics from those in physics. Furthermore, Riemann's philosophical ideas are often in the background of his work on science. 

The aim of this chapter is to give an overview of  Riemann's mathematical results based on physical reasoning or motivated by physics. We also elaborate on the relation with philosophy. While we discuss some of Riemann's philosophical points of view, we review some ideas on the same subjects emitted by Riemann's predecessors, and in particular Greek philosophers, mainly the pre-socratics and Aristotle.

 \bigskip



\end{abstract}

\maketitle

\tableofcontents 

 \section{Introduction} \label{s:intro}

Bernhard Riemann\index{Riemann, Bernhard (1826--1866)} is one of these pre-eminent scientists who considered mathematics, physics and philosophy as a single subject,  whose objective is part of a continuous quest for understanding the world. 
His writings not only are the basis of some of the most fundamental mathematical theories that continue to grow today, but they also effected a profound transformation of our knowledge of nature, in particular through the physical developments to which they gave rise, in mechanics, electromagnetism, heat, electricity, acoustics, and other topics. In Riemann's writings, geometry is at the center of physics, and physical reasoning is part of geometry. His ideas on space and time affected our knowledge in a profound way. They were at the basis of several elaborate theories by mathematicians and physicists, and one can mention here the names of  Hermann Weyl\index{Weyl, Hermann (1885--1955)} and Albert Einstein.\index{Einstein, Albert (1879---1955)} Likewise, Riemann's speculations on the infinitely small and the infinitely large go beyond the mathematical and physical setting, and they had a non-negligible impact on  philosophy.

In the present chapter, we survey some of Riemann's ideas from physics that are contained in his mathematical works. It is not easy to separate Riemann's ideas on physics from those on mathematics. It is also a fact that one cannot consider the fundamental questions that Riemann addressed on physics without mentioning his philosophical background. This is why our survey involves philosophy, besides physics and mathematics. 
 We mention by the way that a certain number of papers and fragments by Riemann on philosophy,  psychology,\index{psychology} metaphysics and gnosiology were collected by Heinrich Weber and published in his edition of Riemann's  Collected Works (p. 507--538). We also mention the name of Gilles Deleuze (1925--1995),\index{Deleuze, Gilles (1925--1995)} a twentieth-century French philosopher who was influenced by Riemann. The name of Deleuze is not commonly known to mathematicians. The relation of his work with Riemann's ideas is highlighted in two chapters of the present volume (see \cite{Jedr} and \cite{Plot}).

As a mathematician, physicist and philosopher,  Riemann belongs to a long tradition of scholars which can be traced back to ancient Greece. 
One of the main outcomes of his Riemann's Habilitation lecture  \emph{\"Uber die Hypothesen, welche der Geometrie zu Grunde liegen} (On the hypotheses that lie at the bases of geometry)\index{Riemann! habilitation lecture}\index{habilitation lecture!Riemann} \cite{Riemann-Ueber} (1854), which we discuss more thoroughly in \S\,\ref{s:vortrag} of the present chapter, is the merging of philosophy, geometry and physics. The fundamental questions that he addresses explicitly in this work, on space,\index{space} form, dimension, magnitude, the infinite and the infinitesimal, the discrete and the continuous are precisely the  questions  that obsessed the Greek philosophers, starting with the Milesians and the Pythagoreans, and passing through Plato,\index{Plato (5th--4th c. B.C.)} Aristotle, Archimedes and several others. One important fact to recall is that the Greeks had a \emph{name} for infinity,\index{infinity} \emph{apeiron}.\index{apeiron@\emph{apeiron}} The name was used by Anaximander,\index{Anaximander (c. 610-- c. 546 B.C)} in the sixth century B.C. There is an extensive literature on the word \emph{apeiron}, whether it denoted an unlimited extent, or a boundless shape, whether it applies to quantity or to shape, etc. The Greeks thoroughly considered the question of infinity, both mathematically and philosophically, and it is often difficult to make the distinction between the two points of view.
A mathematical method of dealing with infinitely small quantities called ``method of exhaustion,"\index{method of exhaustion} which is very close to what we use today in infinitesimal calculus, was developed in the fourth century B.C. by Eudoxus of Cnidus,\index{Eudoxus of Cnidus (4th c. B.C.)} a student of Plato. This method is used by Euclid in the proofs of several propositions of the \emph{Elements}. Dedekind was inspired by this method when he introduced the so-called Dedekind cuts. It is also well known that the philosophical reflections on the infinitely small are not foreign to Leibniz's and Newton's work on the foundations of infinitesimal calculus.

A certain number of these thinkers wondered about the smallest particles of matter, for which they invented a name: atoms, they speculated about their shape and their arrangement and how they fit in an ambient space,\index{space} they meditated on characters of these atoms: cold, waterly, etc.  The thinkers belonging to the ``atomist"\index{atomism} tradition believed that the universe is a mixture of such atoms, that is, uncuttable, or indivisible matter, and void.\index{void} Riemann had his own ideas about matter and void.\index{void} Klein, in his \emph{Development of mathematics in the 19th century} \cite{Klein-development} (p. 235), reporting on some of Riemann's ideas from his \emph{Nachlass} (the collection of manuscripts, notes and correspondence that he left),\index{Riemann's Nachlass} writes:
  \begin{quote}\small
  Riemann thinks of space as being filled with continuous matter \emph{[Stoff]}, which transmits the effect of gravity, light, and electricity. He has throughout the idea of a temporal extension of process. A remark on this topic is found in a personal letter from Gauss to Weber -- with an express request for complete secrecy. And now I again ask, how did these things come to Riemann? It is just mystical influence, which cannot be defined and yet cannot be clearly grasped, of the general atmosphere of a receptive spirit. 
  \end{quote}
 Long before Riemann, the pre-socratics Parmenides\index{Parmenides (6th c. B.C.)} and Zeno\index{Zeno (6th c. B.C.)} (sixth century B.C.), and then  Leucippius\index{Leuceppius (5th c. B.C.)} and Democritus\index{Democritus (c. 460--c. 370 B.C.)} and other scholars of the fifth century B.C., thoroughly gazed at the notions of atom and indivisible matter. Their opinions are reported on by Aristotle,\index{Aristotle (384--322 B.C.)} who made a systematic study of this matter in several texts (\emph{Metaphysics} V, \emph{Physics} V and VI, \emph{Categories} IV, etc.). 
 Other Greek thinkers considered that matter is continuous, rather than atomic, asserting that the atomic structure requires the existence of a void, and claiming that the existence of a void contradicts several laws of physics. They stressed instead the geometric structure of the universe. A theory of chaos,\index{chaos} in the sense of unformed matter arising from the void had also its supporters -- Chaos is an important notion in Greek mythology-- but in general, the Greek philosophers considered that nature is governed by natural laws which they tried to understand. Concerning these thinkers, let us quote Hermann Weyl,\index{Weyl, Hermann (1885--1955)} one of the best representatives of Riemann's tradition of thought, from the beginning of his book \emph{Philosophy of mathematics and natural science}  \cite{Weyl-philo} (p. 3):
\begin{quote}\small
To the Greeks we owe the insight that the structure of space,\index{space} which manifests itself in the relations between spatial configurations and their mutual lawful dependences, is something entirely rational.
\end{quote}
Thus, talking about the origin of Riemann's ideas, we shall often mention his Greek predecessors.

One should also recall that the exceptional rise of Greek science that started in the sixth century B.C.,  in the form of precise questions whose aim was to understand the universe, was accompanied by a profound philosophical reflection on the nature and the goal of sciences, and in particular mathematics. Aristotle,\index{Aristotle (384--322 B.C.)} who is probably the best representative of the Greeks thinkers of the culminating era, in  Book VI of his \emph{Metaphysics} \cite{Aristotle-Metaphysics}, states that among the sciences, three have the status of being \emph{theoretical}: mathematics, physics and theology, the latter, for him, being close to what we now understand as philosophy.\footnote{Cf. \cite{Aristotle-Metaphysics} p. 1619: ``And since natural science, like other sciences, confines itself to one class of being, i.e. to that sort of substance which has the principle of its movement and rest present in itself, evidently it is neither practical nor productive. For the principle of production is in the producer --  it is either reason or art or some capacity, while 
the principle of action is in the doer -- viz. choice, for that which is done and that which is chosen are the same. Therefore, if all thought is either practical or productive or theoretical, natural science must be a theoretical, but it will theorize about such being as admits of being moved, and about substance which in respect of its formula is from the most part not separable from matter.  Now, we must not fail to notice the nature of the essence of its formula, for, without this, inquiry is but idle. [...] 
That natural science, then, is theoretical, is plain from these considerations. Mathematics also, however, is theoretical; but whether its objects are immovable and separable from matter, is not at present clear; it is clear, however, that it considers some mathematical objects \emph{qua} immovable and \emph{qua} separable from matter. But if there is something which is eternal and immovable and separable, clearly the knowledge of it belongs to a theoretical science, -- not, however, to natural science (for natural science deals with certain movable things) nor to mathematics, but to a science prior to both."} Let us note right away that these are precisely the three branches of knowledge that constitute the background of Riemann (who, by the way, was also trained in theology). We also note that although Pythagoras is supposed to have coined the term $\phi \iota  \lambda o \sigma o  \phi \acute{\iota} \alpha $,\index{philosophy@$\phi \iota  \lambda o \sigma o \acute{\iota} \phi \alpha $} \emph{communis opinio} now seems that its current meaning (striving for knowledge) goes back to Plato.\index{Plato (5th--4th c. B.C.)}

 In the same work and in others, Aristotle\index{Aristotle (384--322 B.C.)} discusses at length the role of each of these three sciences and the relations among them. He also addresses thoroughly the question of whether mathematics has a purely ideal character or whether it  reflects the real world. Such interrogations lead directly to the most fundamental questions that Riemann addressed in his Habilitationsvortrag\index{Riemann! habilitation lecture}\index{habilitation lecture!Riemann} and in his other writings. We shall say more on the lineage of Riemann's ideas to Greek philosophy in \S \ref{s:vortrag}  where we discuss this work.

Riemann's interest in physics was constant during his lifetime. Since his early twenties, he tried to develop a theory that would unify electricity, magnetism, light and gravitation -- the same quest that Poincar\'e, Lorentz and Einstein had after him, culminating in the theory of general relativity.  One of Riemann's  manuscripts, \emph{Ein Beitrag zur Electrodynamik}  (A contribution to electrodynamics) \cite{Riemann-Heat},  whose subject is electrodynamics\index{electrodynamics} and which is related to his search for the unification of the various  forces of nature, published posthumously, is discussed in the chapter \cite{Goenner} by Hubert Goenner in the present volume. In this paper, Riemann develops a theory of electromagnetism which is based on the assumption that electric current travels at the velocity of light. Furthermore, he considers that the differential equation that describes the propagation electric force is the same as that of heat and light propagation.
Goenner, in that paper, mentions the works of other physicists of the same period, including Maxwell,\index{Maxwell, James Clerk (1831--1879)} Lorenz,  Helmholtz,\index{Helmholtz@von Helmholtz, Hermann Ludwig Ferdinand (1821--1894)} Carl Neumann and Franz Neumann.  Being himself a physicist, Goenner writes:  
\begin{quote}\small
Surprisingly, within the then reigning view of electromagnetism as a
particle theory, we can note a relativistic input, made by the famous
mathematician Bernhard Riemann. His introduction of the
retarded scalar potential into theoretical electrodynamics is still
valid, but remains unknown to the overwhelming majority of today's
theoretical physicists.
\end{quote}
Several other commentaries on Riemann's manuscript exist, and some of them are mentioned in the bibliography of \cite{Goenner}. Enrico Betti,\index{Betti, Enrico (1823--1892)} who was Riemann's friend and who translated into Italian several of his works and wrote commentaries on them, had already commented on that paper in 1868, see \cite{Betti-elettrodinamica}. 

It is known that, as a student, in G\"ottingen and  Berlin, Riemann
attended more courses and seminars on (theoretical and experimental) physics than on mathematics. One may also mention here an essay in Riemann's \emph{Nachlass}, entitled \emph{Gravitation und Licht} (Gravitation and light) \cite{Riemann-Gesammelte} p. 532--538, whose subject is the theoretical connection between gravitation and light.  Betti\index{Betti, Enrico (1823--1892)} also wrote a paper, entitled \emph{Sopra una estensione dei principii generali della dinamica} (On the extension of the general principles of dynamics) \cite{Betti-sopra}, in which he announces several results which are based on ideas contained in Riemann's lectures \emph{Schwere, Electricit\"at und Magnetismus}, edited by the latter's student Karl Hattendorf (Hannover, 1880)\index{Hattendorf, Karl (1834--1882)}  \cite{Riemann-Schwere}. Riemann establishes in these lectures necessary and sufficient conditions under which Hamilton's principle on the motion of a free system subject to time-independent forces that depend on the position and the motion of the system is satisfied. Chapter VII of Picard's\index{Picard, Charles \'Emile (1856--1941)} famous \emph{Trait\'e d'analyse} \cite{Picard-traite} contains a chapter called \emph{Attraction and potential}. The author declares there (p. 167) that he uses a transformation from Riemann's posthumous memoir  \emph{Schwere, Electricit\"at und Magnetismus}. Finally, we note that Maxwell\index{Maxwell, James Clerk (1831--1879)} discussed Riemann's theory of electrodynamics in his \emph{Note on the electromagnetic theory of light}, an appendix to his paper \cite{Maxwell1868}. 

 In talking about Riemann's background in physics, we take this opportunity to recall a few facts about Riemann's studies. 
 
In a letter to his father on April 30, 1845, while he was still in high-school, Riemann informs the latter that he starts being more and more attracted by mathematics.  He also tells him  in the same letter that he plans to enroll the University of G\"ottingen to study theology, but that in reality he must decide for himself what he shall do, since otherwise he will not bring anything good to a subject. Cf. \cite{Riemann-Letters}. Riemann entered the University of G\"ottingen in 1846, as a student in theology. He stayed there for one year and then moved to the University of Berlin where he spent two years, attending lectures by Jacobi\index{Jacobi, Carl Gustav Jacob (1804--1851)} and Dirichlet. In a letter to his father, dated May 30, 1849, Riemann writes (\cite{Riemann-Letters}):\footnote{We are using the translation by Gallagher and  Weissbach.} ``I had come just in time for the lectures of Dirichlet and Jacobi. Jacobi has just begun a series of lectures in which he lead off once again with the entire system of the theory of elliptical functions in the most advanced, but elementary way." Jacobi was highly interested in mechanics, and it would not be surprising if his interest in elliptic functions was motivated by their applications to mechanics. In another letter, written to his brother, dated November 29, 1847, Riemann writes: 
\begin{quote}\small
When I arrived, I found to my great joy that Jacobi, who had announced no course in the catalog, had changed his mind. He plans to lecture on mechanics. I would, if possible, stay here for another semester just to attend it. Nothing could be more satisfying to me than this. [...] The next day I went to see Jacobi in order to enroll in his course. He was very polite and friendly, because in his previous lecture he had dealt with a  subject related to the problem I had just solved, I brought it up and told him of my work. He said if it was a nice job, he would send it to Crelle's Journal as soon as possible. Unfortunately my time will be somewhat tight for writing it up. Also I dont know whether the complete solution of the problem will take yet more time.\footnote{It is not clear to the author of the present article what the work that Riemann is talking about is. It might be that there is a mistake in the date of the letter. Nevertheless, the content is interesting for us here regardless of the date.}
\end{quote}

Dirichlet's lectures in Berlin, at that epoch, were centered mostly on theoretical physics (partial differential equations). It is 
from these lectures that Riemann became familiar with potential theory,\index{potential theory} a topic which was about to play an important role in his later work. Klein, in his \emph{Development of mathematics in the 19th century} (\cite{Klein-development}, p. 234--235), writes:
\begin{quote}\small
Dirichlet loved to make things clear to himself in an intuitive substrate; along with this he would give acute, logical analyses of foundational questions and would avoid long computations as much as possible. His manner suited Riemann, who adopted it and worked according to Dirichlet's methods.
\end{quote}
Riemann's admiration for Dirichlet is expressed at several places of his writings, for instance in the third section of the historical part of his habilitation dissertation on trigonometric series which we shall analyse in \S\ref{s:trigo} below.\footnote{It is also true that Riemann, with his extreme sensibility, was at some point disappointed of Dirichlet being less amicable with him. In a letter to his brother, dated April 25, 1857, he writes: [...] Also Dirichlet appeared, if still very polite, yet not so well-disposed towards me as before. This also was agony for me.}

During his stay at Berlin, Riemann also attended lectures by Dove on optics, which he found very interesting, and by Enke on astronomy (letter without date, quoted in \cite{Riemann-Letters}). About the latter, Riemann says that ``his presentation for the most part is rather dry and boring, however the time that we spend at the observatory once a week, from 6pm to 8pm, is useful and instructive."

 After the two years spent in Berlin (1847--1849), Riemann returned to G\"ottingen, where  he attended the lectures and seminars of the newly hired physics professor Wilhelm Eduard Weber (1804-1891),\footnote{This was the second time that Weber was hired in G\"ottingen. The first time was in 1831, at Gauss's recommendation (Weber was 27). Weber developed a theory of electromagnetism which was eventually superseded by Maxwell's.\index{Maxwell, James Clerk (1831--1879)} Weber and Gauss published joint results and they constructed the first electromagnetic telegraph (1833), which operated between the astronomical observatory and the physics laboratory of the University of G\"ottingen (the locations were 3 km apart). In 1837, as the result of a repression, led by the new King of Hannover Ernest Augustus (who reigned between 1837 and 1851) and caused by political events, Weber, together with six other leading professors (including the two brothers Grimm), was dismissed from his position at the University of G\"ottingen. He came back to this university in 1849 and served intermittently as the administrator of the Astronomical Observatory. The position of director had been occupied by Gauss since its foundation in 1816. Gauss was more than seventy at the time Weber returned to G\"ottingen.}\index{Weber,Wilhelm Eduard (1804--1891)} who was also Gauss's collaborator and close friend. Klein writes, in his \emph{Development of mathematics in the 19th century}, (\cite{Klein-development}, p. 235): ``In Weber, Riemann found a patron and a fatherly friend. Weber recognized Riemann's genius and drew the shy student to him. [...] Riemann's interest in the mathematical treatment of nature was awakened by Weber, and Riemann was strongly influenced by Weber's questions."  

 From Riemann's posthumous papers, we read:\footnote{Cf. Bernhard Riemann's \emph{Gesammelte mathematische Werke und wissenschaftlicher Nachlass}, ed. H. Weber and R. Dedekind, \cite{Riemann-Gesammelte} 2nd edition, Leipzig, 1892, p. 507. The translation of this passage is borrowed from the English translation of Klein's \emph{Development of mathematics in the 19th century}, (\cite{Klein-development}, p. 233).}  
\begin{quote}\small 
My main work concerns a new conception of the known works of nature -- their expression by means of other basic concepts -- whereby it became possible to use the experimental data on the reciprocal actions between heat, light, magnetism and electricity to investigate their connections with each other. I was led to this mainly through studying the works of Newton,\index{Newton, Isaac (1643--1727)} Euler,\index{Euler, Leonhard (1707--1783)} and -- from another aspect -- Herbart.\index{Herbart, Johann Friedrich (1776--1841)}
\end{quote}

The name of Herbart\index{Herbart, Johann Friedrich (1776--1841)}, which appears at the end of this quote, will be mentioned several times in the present chapter, and it is perhaps useful to say right away a few words on him.

 Johann Friedrich Herbart (1776--1841)\index{Herbart, Johann Friedrich (1776--1841)}  started his studies in philosophy in Iena under Fichte\index{Fichte,  Johann Gottlieb (1762--1814)}  but he soon disagreed with his ideas and went to G\"ottingen  where he received his doctorate and habilitation, and after that he taught there pedagogy and philosophy. In 1809, he was offered the chair formerly held by Kant\index{Kant, Immanuel (1724--1804)} in K\"onigsberg.  His philosophy relies on Leibniz's\index{Leibniz, Gottfried Wilhelm (1646--1716)} theory of monads. Herbart\index{Herbart, Johann Friedrich (1776--1841)} was conservative and anti-democratic. He was an advocate of the view that the state higher officials should be appointed among those who have a strong cultural education. He wrote in 1824 a treatise entitled \emph{Psychology as a science newly founded on experience, metaphysics and mathematics} \cite{Herbart}. In his research on psychology,\index{psychology} he made use of infinitesimal calculus, and he was probably the first to do this. For him, psychology is a science which is based at the same time on experimentation, mathematics and metaphysics, and he made a parallel between this new field and the field of physics in the way Newton\index{Newton, Isaac (1643--1727)} conceived it. Sygmund Freud\index{Freud, Sygmund (1856--1939)} was profoundly influenced by Herbart.\footnote{From the article on Herbard in the \emph{Freud encyclopedia: Theory, therapy and culture} (\cite{Freud} p. 254), we read: ``The ghost of the philosopher Johann Friedrich Herbart\index{Herbart, Johann Friedrich (1776--1841)} hovers over all of Freud's works, an inseparable albeit unacknowledged presence. Herbart, the successor of Kant\index{Kant, Immanuel (1724--1804)} in K\"onigsberg, arguably exercised a more profound, more persuasive influence on Freud\index{Freud, Sygmund (1856--1939)} than either Schopenhauer\index{Schopenhauer, Arthur (1788--1860)} or Nietzsche,\index{Nietzsche, Friedrich (1844--1900)} whom many scholars regard as sources for some of his major concepts. From Herbart,\index{Herbart, Johann Friedrich (1776--1841)} Freud derived such ideas as the mental activity can be conscious, preconscious, or unconscious, that unconscious mental activity is a continuous determinant of conscious activity, and that the present is unceasingly shaped by the past, whether remembered or forgotten. From Herbart, he also borrowed some essentials of his model, the idea of conflicting conscious and unconscious psychic forces, the censorship-exercising ego, the threshold of consciousness, `resistance,' `repression' and much else. [...] It was Herbart's\index{Herbart, Johann Friedrich (1776--1841)} ambition to contribute to the establishment of `a research of the mind which will be the equal of natural science, insofar as this science everywhere presupposes the absolutely regular connection between appearances.' He compared the situation in psychology\index{psychology} with that of astronomy: in the pre-Copernician era,\index{Copernican revolution} the motions of the planets had seemed irrational; every so often these heavenly bodies inexplicably seemed to reverse their course; for this reason they were known as the `wanderers.' These peculiar paths, however, were recognized as entirely lawful as soon as the heliocentric theory was introduced. The hypothesis of unconscious thought\index{unconscious thought} performed the same service to the mind, Herbart\index{Herbart, Johann Friedrich (1776--1841)} maintained." We encourage the interested reader to go through this entire article by R. Sand.} 
 Herbart returned to G\"ottingen in 1833 where he taught philosophy and pedagogy until his death. Riemann was 15 and it is unlikely that he followed any course of Herbart.\index{Herbart, Johann Friedrich (1776--1841)}
 Erhard Scholz, in his paper \cite{Scholz1982}, reports on several sets of notes written by Riemann and preserved in the Riemann archives in G\"ottingen, which concern the philosophy of Herbart. These notes show that Riemann was indeed influenced by the philosopher, for what concerns epistemology and the philosophy of science, and in particular for his ideas on space.\index{space} After analyzing some of Riemann's fragments, Scholz writes: ``Riemann's views on mathematics seem to have been deepened and clarified by his extensive studies of Herbart's philosophy. Moreover, without this orientation, Riemann might have never formulated his profound and innovative concept of manifold. This represents an indirect but nevertheless effective influence of Herbart on Riemann's mathematical and (in particular) his geometrical thinking." The last paragraph we quoted from Riemann's  \emph{Gesammelte mathematische Werke} continues as follows:
 \begin{quote}\small
 As for the latter, I could almost completely agree with Herbart's earliest investigations, whose results are giving in his graduating and qualifying theses \emph{[Promotions --  und Habilitationsthesen]} (of October 22 and 23, 1802), but I had to veer away from the later course of his speculations at an essential point, thus determining a difference with respect to his philosophy of nature and those propositions of psychology which concerns its connections with the philosophy of nature.
 \end{quote}

When Riemann became, in 1854, \emph{Privatdozent} at the University of G\"ottingen, the subject of the first lessons he gave was  differential equations with applications to physics. 
These lectures became well known even outside Germany. 
In a letter to Hou\"el\index{Hou\"el, Jules (1823--1886)}  sent in 1869 (see \cite{DH} p. 90), Darboux\index{Darboux, Gaston (1842--1917)} writes:\footnote{In the present chapter, the translations from the French are mine, except if the contrary is indicated.} 
\begin{quote}\small
I wonder if you know a volume by Riemann on mathematical physics entitled \emph{On partial differential equations}. I was very much interested in this small volume, it is clear and could be put in the hands of the students of our universities. I think that you will appreciate it; if you are a little bit concerned with mathematical physics it will be of interest to you. Above all, hydrodynamics seems to me very well treated.\footnote{[Je ne sais si vous connaissez un volume de Riemann sur la physique math\'ematique intitul\'e sur les \'equations aux d\'eriv\'ees partielles. Ce petit volume m'a beaucoup int\'eress\'e, il est clair et pourrait \^etre mis avec avantage entre les mains des auditeurs de nos facult\'es. Je crois que vous en serez content ; si vous vous occupez un peu de physique math\'ematique il vous int\'eressera, l'hydrodynamique surtout m'y a paru tr\`es bien trait\'ee.]}
\end{quote}

Three editions of Riemann's notes from his lectures on differential equations applied to physics  appeared in print, the last one in 1882, edited by his student  Karl Hattendorf.\index{Hattendorf, Karl (1834--1882)} A work in two volumes entitled 
\emph{Die partiellen Differential-Gleichungen der mathematischen Physik}  (The partial differential equations of mathematical physics) \cite{Riemann-Weber} appeared in 1912, written by Heinrich Weber.\footnote{Heinrich Weber (1842--1913)\index{Weber, Heinrich (1842--1913)}  taught mathematics at the University of Strasbourg -- a city which at that time belonged to Germany -- from 1895 till 1913. He is the co-editor, with Dedekind, of Riemann's collected works. Dedekind\index{Dedekind, Richard (1831--1916)} followed Riemann's lectures in 1855--1856, and they became friends. Regarding the relation between the two mathematicians, let us mention the following. In 1858, a position of professor of mathematics was open at the Zurich Eidgen\"ossische Technische Hochschule. Both Riemann and Dedekind applied, and Dedekind was preferred, probably because he had more experience in teaching elementary courses. Indeed, a Swiss delegation visited G\"ottingen to examine the candidates, and considered that Riemann was ``too introverted to teach future engineers." After he left for Zurich, Dedekind remained faithful to Riemann. Klein, in his Lectures on the history of nineteenth century mathematics, characterized Dedekind as a major representative of the Riemann tradition. At Riemann's death in 1866, Dedekind was given the heavy load of editing Riemann's works. This is where he asked the help of his friend Heinrich Weber. The first edition of the Collected works was published in 1876, and it included a biography of Riemann written by Dedekind. Before that, Dedekind had edited, in 1868, the two habilitation works of Riemann, \emph{\"Uber die Darstellbarkeit einer Function durch eine trigonometrische Reihe} (On the representability of a function by a trigonometric series)\index{trigonometric series} \cite{Riemann-Trigo} and \emph{\"Uber die Hypothesen, welche der Geometrie zu Grunde liegen} \cite{Riemann-Ueber} which we already mentioned. Klein says in his lectures on the history of nineteenth century mathematics \cite{Klein-development} that at Riemann's death, the latter's heirs entrusted Dedekind with the edition of the \emph{Nachlass},  that Dedekind started working on that and wrote illuminating comments, but that he was not able to continue that work alone. In 1871 he asked the help of Clebsch,\index{Clebsch, Alfred (1833--1872)} but the latter died soon after (in 1872). He then asked Weber to help him completing the work. In 1882, Dedekind\index{Dedekind, Richard (1831--1916)}  and Weber\index{Weber, Heinrich (1842--1913)} published a paper entitled \emph{Theorie der algebraischen Functionen einer Veranderlichen} (The theory of algebraic functions in one variable) \cite{DW}  in which they developed in a more accessible manner Riemann's difficult ideas on the subject. An analysis of this paper, including a report on its central place in the history of mathematics, is contained in Dieudonn\'e \cite{Dieudonne1985} p. 29--35.
Finally let us mention that Weber is also a co-editor of the famous Klein \emph{Encyklop\"adie der Mathematischen Wissenschaften}.}   This work is considered as a revised version of Riemann's lectures on this subject. The applications to physics include heat conduction, elasticity and hydrodynamics.  The work was used for many years as a textbook in various universities.\footnote{See \cite{Ames} for a review of this book.} In a biography of Neugebauer \cite{Rowe} (p. 16), the author reports that in the 1920s, Courant was using these books in his teaching:
 \begin{quote}\small
 Courant\index{Courant, Richard (1888--1972)} was hardly a brilliant lecturer, but he did have the ability to spark interest in students by escorting them on the frontiers of research in analysis. In his course on partial differential equations, he stressed two sharply opposed types of literature: works that expound general theory, on the one hand, and those that pursue special problems and methods, on the other [...] For the second type of literature, Courant's background and personal preferences came to the fore. Here he recommended the \emph{Ausarbeitung} of Riemann's lectures on partial differential equations prepared by Karl Hattendorf\index{Hattendorf, Karl (1834--1882)} together with Heinrich Weber's subsequent volume on Riemann's theory of PDEs in mathematical physics.
\end{quote}

Concerning physics in Riemann's work, let us also quote Klein\index{Klein, Felix (1849--1925)} from his address delivered at the 1894 general session of the \emph{Versammlung Deutscher Naturforscher und \"Artzte} (Meeting of the German naturalists and physicians; Vienna, September 27, 1894), cf. \cite{Klein1923}:
 \begin{quote}\small
I must mention, first of all, that Riemann devoted much time and thought to physical considerations. Grown up under the tradition which is represented by the combinations of the names of Gauss and Wilhelm Weber, influenced on the other hand by Herbart's\index{Herbart, Johann Friedrich (1776--1841)} philosophy, he endeavored again and again to find a general mathematical formulation for the laws underlying all natural phenomena. [...] The point to which I wish to call your attention is that \emph{these physical views are the mainspring of Riemann's purely mathematical investigations}.
\end{quote}

\medskip

 Riemann wrote papers on physics (a few, but may be as many as his papers on mathematics), and they will only be briefly mentioned in the present chapter.  
Our main subject is \emph{not} Riemann's work on physics, but his ideas concerning physics that are present in his mathematical papers. In these ideas, we shall comment on his motivation. 
 Jeremy Gray, the author of the chapter of the present volume entitled \emph{Riemann on geometry, physics, and philosophy} \cite{Gray-R}, writes in that chapter that ``Riemann belongs to a list of brilliant mathematicians whose lasting contributions are more in mathematics than physics, contrary to their hope."
We shall mention Riemann's impact of some of his mathematical work on the later development of physics. At the same time, we shall give an overview of some important mathematical works of Riemann, in particular on the following topics:
\begin{enumerate}
\item \emph{Riemann surfaces and functions of a complex variable}. Riemann approached complex analysis from the point of view of potential theory, that is, based on the theory of the Laplace equation\index{Laplace equation}\index{equation!Laplace} 
\[\frac{\partial^2 u}{\partial x^2} + \frac{\partial^2 u}{\partial y^2}=0
.\]
Here, the function $u=u(x,y)$ represents the potential function that gives rise to a  streaming for an incompressible flow contained between two planes parallel to the $x,y$ plane.  (the flow may be an electric field, in which case $u$ is the electrostatic potential). The Laplace equation expresses the fact that there is as much fluid that flows into an element of area per unit time than fluid that flows out.
The bases of the theory of Riemann surfaces are contained in Riemann's doctoral dissertation\index{doctoral dissertation!Riemann}\index{Riemann!doctoral dissertation} \cite{Riemann-Grundlagen} and his paper on Abelian functions \cite{Riemann-Abelian}. We shall review the role of physical reasoning in these works. Conversely, Riemann made the theory of functions of a complex variable, based on his approach using partial differential equations in the real domain, a basic tool in mathematical physics. We shall mention below some of the tremendous impact of the theory of Riemann surfaces in modern theoretical physics.

\item \emph{Trigonometric series}:\index{trigonometric series} Riemann's work on trigonometric\index{trigonometric series} series is contained in his habilitation dissertation (\emph{Habilitationsschrift})\index{Riemann! habilitation text}\index{habilitation text!Riemann} \cite{Riemann-Trigo}. His motivation, as Riemann himself writes, comes from the theory of sound. The origin of the questions he tackled lies in seventeenth-century physics and mathematics, and they led then to a harsh debate that involved several scientists including Euler,\index{Euler, Leonhard (1707--1783)} Lagrange,\index{Lagrange, Joseph-Louis (1736--1813)}  d'Alembert\index{Alembert@d'Alembert, Jean le Rond (1717--1783)} and Daniel Bernoulli (to mention only the most famous ones). From the mathematical point of view, the main issue was the nature of the functions that were admitted as solutions of the wave equation.\index{wave equation}\index{equation!wave} Riemann eventually concluded the debate, showing the generality of  those functions that have to be included as solutions of these equations. In the same paper, Riemann laid the foundations of what became known later on as the theory of the Riemann integral. This came from his effort to clarify the nature of the coefficients of a trigonometric series\index{trigonometric series} associated with a function. These coefficients are indeed given in the form of integrals.

\item \emph{Riemannian geometry}. This is contained in Riemann's habilitation lecture \cite{Riemann-Ueber}  and his later paper, the \emph{Commentatio} \cite{Commentatio}.\index{Riemann! habilitation text}\index{habilitation text!Riemann} In the development of this theory, Riemann was motivated in part by physics, and in part by philosophy. In his habilitation lecture, Riemann's bond of filiation with Greek philosophy, and in particular with Aristotle,\index{Aristotle (384--322 B.C.)} is clear. We shall comment on this and we shall also recall the huge impact of these two works of Riemann on the later physical theories. 

\item \emph{Other works}. In the last section of this chapter, we shall analyze more briefly some other papers of Riemann related to our subject. 

\end{enumerate}

To close this introduction, we mention that the fact that Riemann, in his mathematical work, was motivated by physics was also common to other mathematicians of the eighteenth and the nineteenth centuries.  One may recall that Gauss, who was Riemann's mentor, considered himself more as a physicist than a mathematician. We refer to \cite{Garland} for a review of Gauss's contribution to geomagnetism. Gauss\index{Gauss, Carl Friedrich (1777--1855)} was in charge of the practical task of surveying geodetically the German kingdom of Hannover. In the preface of his paper \cite{Gauss-Copenhagen} which we already mentioned, published in 1825, about the same time he wrote his famous \emph{Disquisitiones generales circa superficies curvas} (General investigation of curved surfaces) (1825 and 1827) \cite{Gauss-English}, Gauss writes that his aim is only to construct geographical maps and to study the general principles of geodesy\index{geodesy} for the task of land surveying. Surveying  the kingdom of Hannover took  nearly two decades to be completed. It led Gauss gradually to the investigation of triangulations,  to the use of the method of least squares in geodesy,\footnote{Gauss first published his method of least squares in an important treatise in two volumes calculating the orbits of celestial bodies in 1809 \cite{Gauss1809}, but in that work he claims that he knew the method since 1795. This led to a priority controversy between Gauss\index{Gauss, Carl Friedrich (1777--1855)} and Legendre,\index{Legendre, Adrien-Marie (1752--1833)} who published the first account of that method in 1805 \cite{Legendre}.} and then to his \emph{Disquisitiones generales circa superficies curvas}. In the latter, we can read, for instance, in 
\S 27 (p. 43 of the English translation \cite{Gauss-English}):
``Thus, e.g., in the greatest of the triangles which we have measured in recent years, namely that between the points Hohenhagen, Brocken, Inselberg, where the excess of the sum of the angles was 14."85348, the calculation gave the following reductions to be applied to angles: 
Hohehagen: 4."95113; 
Brocken: 4."95104;
Inselberg: 4."95131."

It is also interesting to know that Jacobi,\index{Jacobi, Carl Gustav Jacob (1804--1851)} after Gauss, studied similar problems of geodesy,\index{geodesy} using elliptic functions.\index{function!elliptic}\index{elliptic function}  In a paper entitled \emph{Solution nouvelle d'un problème fondamental de géodésie} (A new solution of a fundamental problem in geodesy) \cite{Jacobi-geodesie}, he considers, on an ellipsoid having the shape of the earth, a geodesic arc whose length, the latitude of its origin and its azimuth angle at that point are known. The question is then to find the latitude, the azimuth angle of the extremity of this arc, as well as the difference in longitudes between the origin and the extremity. He then declares: ``The problem of which I just gave a new solution has been recently the subject of a particular care from Mr. Gauss, who treated it in various memoirs and gave different solutions of it."\footnote{Le problème dont je viens de donner une solution nouvelle a été dans ces derniers temps l'objet de soins particuliers de la part de M. Gauss, qui en a traité dans différents mémoires et en a donné plusieurs solutions.}

Riemann was profoundly influenced by Gauss. We emphasize this fact because it is written here and there that Riemann did not learn a lot from Gauss, since when Riemann started his studies, Gauss was already old, and that  in any case, Gauss was never interested in teaching. Klein writes in his  \emph{Development of mathematics in the 19th century}  (\cite{Klein-development}, p. 234 of the English translation): 
\begin{quote}\small
Gauss taught unwillingly, had little interest in most of his auditors, and was otherwise quite inaccessible. Nevertheless, we call Riemann a pupil of Gauss; indeed he is Gauss's only true pupil, entering into his inner ideas, as we now are coming to see in outline from the \emph{Nachlass}.
\end{quote}

 Before Gauss, Euler,\index{Euler, Leonhard (1707--1783)} whose work was also a source of inspiration for Riemann,  was likewise thoroughly involved in physics. His work on partial differential equations was motivated by problems from physics. In fact, Euler tried to systematically reduce every problem in physics to the study of a differential equation.  Euler was also very much involved in acoustics.  The initial attraction by Euler to number theory arose in his work on music theory; cf. \cite{Euler-Hermann} where this question is thoroughly discussed. Later in this chapter, we shall have the occasion to talk about Euler's and Riemann's works related to acoustics. 
The influence of Euler and Gauss on Riemann is thoroughly reviewed in Chapter 1 of the present volume \cite{Papa-Euler}.

\section{Function theory and Riemann surfaces}\label{s:RS}

Riemann's work on the theory of functions of a complex variable is developed in his two memoirs  \emph{Grundlagen f\"ur eine allgemeine Theorie der Functionen einer ver\"anderlichen complexen Gr\"osse} (Foundations of a general theory of functions of a variable complex magnitude) \cite{Riemann-Grundlagen} (1851) and \emph{Theorie der Abel'schen Functionen} (Theory of Abelian functions) \cite{Riemann-Abelian} (1857). The first of these memoirs is Riemann's doctoral dissertation.\index{doctoral dissertation!Riemann}\index{Riemann!doctoral dissertation}
The text of this dissertation was submitted to the University of G\"ottingen on November 14, 1851, and the defense took place on December 16 of the same year. 
Several ideas introduced in these two papers are further developed in subsequent works of Riemann. We mention in particular a fragment on the theory of Abelian functions\index{Abelian function}  published posthumously and which is part of Riemann's collected papers editions \cite{Riemann-Gesammelte} and \cite{Riemann--French}. There are also other works of Riemann that involve in an essential way functions of a complex variable; a famous example is his extension of the real zeta function\footnote{The real zeta function was already considered by Euler; cf. in particular his papers \cite{E72} and \cite{E352}. Note that Euler did not use the notation $\zeta$. In his book, \emph{Introductio in analysin infinitorum} (Introduction to the analysis of the infinite) \cite{Euler-Int-b}, where he considers this function for integer values of the variable, he denotes it by $P$.} to the complex domain, which turned out to be a huge step in the study of this function. Finally, we mention that there are lecture notes of Riemann on functions of a complex variable available at G\"ottingen's library, and there is an outline of these lectures in Narasimhan's article \cite{RN}. It is not our intention here to comment on Riemann's  fundamental work on functions of a complex variable and its importance for later mathematics; we shall only concentrate on its relation to physics. However, we start with a few comments on the theory of functions of a complex variable, before Riemann started working on it, because this will help including Riemann's work in its proper context.
  
The notion of function of a complex variable can be traced back to the beginning of the notion of function, which, in the form which is familiar to us today, it is usually attributed to Johann Bernoulli\index{Bernoulli, Johann (1667--1748)} and\index{Euler, Leonhard (1707--1783)} Euler.\footnote{One should emphasize that the seventeenth-century infinitesimal calculus of Leibniz\index{Leibniz, Gottfried Wilhelm (1646--1716)} and Newton,\index{Newton, Isaac (1643--1727)} which was developed before Euler, dealt with curves, and not with functions.} A precise definition of a function, based on a careful description of the notion of variable, is contained in Euler's treatise \emph{Introductio in analysin infinitorum} \cite{Euler-Int-b} (1748).  This book was published one year after the appearance of the famous memoir of d'Alembert\index{Alembert@d'Alembert, Jean le Rond (1717--1783)} \cite{Al_recherches} in which the latter gave the wave equation.\index{wave equation}\index{equation!wave} We mention this fact because the main mathematical question that was motivated by d'Alembert's\index{Alembert@d'Alembert, Jean le Rond (1717--1783)} memoir turned out to be the question of the nature of functions that are solutions of the wave equation.\index{wave equation}\index{equation!wave} Hence, the general question was addressed: \emph{What is a function?}
Furthermore, this memoir of d'Alembert\index{Alembert@d'Alembert, Jean le Rond (1717--1783)} was the original motivation for the study of trigonometric series,\index{trigonometric series} which was the subject of Riemann's \emph{Habilitationsschrift}\index{Riemann! habilitation text}\index{habilitation text!Riemann} which we discuss in \S \ref{s:trigo}. The \emph{Introductio} consists of two volumes. In Chapter 1 of the first volume,  after he defines functions, Euler\index{Euler, Leonhard (1707--1783)} writes: ``[...] Even zero and complex numbers are not excluded from the signification of a variable quantity."  Thus, complex functions were  considered by Euler from the very outset of his work on general functions.  We shall come back in \S \ref{s:trigo} to Euler's definition of a function.   Two years after his \emph{Introductio}, Euler published his famous memoir \emph{Sur la vibration des cordes} (On string vibration) \cite{Euler_vib-corde} which we shall also discuss below.

After Euler, one has to mention Cauchy,\index{Cauchy, Augustin-Louis (1789--1857)} who made a thorough and profound contribution to the theory of analytic functions of a complex variable, during the three decades that preceded Riemann's work on the subject. In a series of \emph{Comptes Rendus} Notes and in other publications, including his \emph{Cours d'analyse de l'\'Ecole Royale Polytechnique} (A course of analysis of the \'Ecole Polytechnique) \cite{Cauchy-cours} (1821) and his \emph{M\'emoire sur les int\'egrales d\'efinies prises entre des limites imaginaires} (Memoir on the definite integrals taken between two imaginary limits) \cite{Cauchy-memoire} (1825), Cauchy\index{Cauchy, Augustin-Louis (1789--1857)}  introduced several fundamental notions, such as the disc of convergence of a power series and path integrals between two points in the complex plane with the study of the dependence on the path.\footnote{One should mention that the idea of integration along paths was present  in the works of Gauss \cite{Gauss-Uber}
 and Poisson \cite{Poisson-Suite}. They both considered line integrals in the complex plane and they noticed that these integrals depend on the choice of a path.} He dealt with functions which may take the value infinity at some points, and he invented the calculus of residues and the characterization of complex analyticity by the partial differential equations satisfied by the real and imaginary parts of the function, which were called later the Cauchy--Riemann equations.\index{Cauchy--Riemann equations}\index{equations!Cauchy--Riemann}

Besides the work of Cauchy, we mention that of his student Puiseux who further developed some of his master's ideas and brought new ones, essentially in two papers \cite{Puiseux} and \cite{Puiseux1}.  In the 177-page paper  \cite{Puiseux}, Puiseux uses the methods introduced by Cauchy on path integration in the study of the problem of uniformization\index{uniformization}\index{problem!uniformization} of an algebraic function $u(z)$. This is a function defined implicitly by an equation of the form $P(u,z)=0$ where $P$ is a two-variable polynomial. The uniformization problem, in this setting, is to get around the fact that such a function $u$ is multi-valued and to make it univalued (uniform). In doing this,
Puiseux also developed the theory of functions of a complex variable which are of the form $\int udz$, where $u$ is above.  He highlighted the role of the critical points of the function $u$ in this line integral, and the fact that integrating along the loops that contain one such point one gets different values for the function. Using this fact, he gave an explanation for the periodicity of the complex circular functions, of elliptic functions, and of the functions defined by integrals introduced by Jacobi. He showed that for a given $z$, the various solutions $u(z)$ of the equation $f(u,z)=0$ constitute a certain number of ``circular systems," and he gave a method to collect them into groups. In doing this, he developed a geometric Galois theory, discussing the ``substitutions" which act on the solutions of the algebraic equation. He also gave a method to find expressions for these solutions as power series with fractional exponents. There is a profound relation between the results of Puiseux on algebraic functions and Riemann surfaces. This is also surveyed in the chapter \cite{Papa-Puisex} in the present volume.

Since we talked about elliptic functions,\index{function!elliptic} whose study was one of Riemann's main subjects of interest, let us mention that these functions were also used in physics, and that this was certainly one of the reasons why Riemann was interested in them. Already Euler, in his numerous memoirs on elliptic integrals, studied their applications to the oscillations of the pendulum with large amplitudes, to the measurement of the earth, and to the three-body problem. 
In the preface of the treatise  \emph{Th\'eorie des fonctions doublement périodiques et, en particulier, des fonctions elliptiques} (Theory of doubly periodic functions, and in particular elliptic functions) by Briot and Bouquet which we review in another chapter of the present volume \cite{Papa-Riemann3} in relation with Riemann's work, the authors write: 
\begin{quote}\small
One encounters frequently elliptic functions in questions of geometry, mechanics, or mathematical physics. We quote, as examples, the ordinary pendulum, the conical pendulum, the ellipsoid attraction, the motion of a solid body around a fixed point, etc. Mr. Lamé published last year a very interesting work, where he shows that that elliptic functions enter into questions relative to heat distribution  and of isothermal surfaces.\footnote{On rencontre fréquemment les fonctions elliptiques dans les questions de géométrie, de mécanique ou de physique mathématique. Nous citerons, comme exemples, le pendule ordinaire, le pendule conique, l'attraction des ellipso\"\i des, le mouvement d'un corps solide autour d'un point fixe, etc. M. Lamé a publié l'année dernière un ouvrage très-intéressant, où il montre que les fonctions elliptiques s'introduisent dans les questions relatives à la distribution de la chaleur et aux surfaces isothermes.}
\end{quote}
Regarding the same subject, we note that the second volume of Halphen's \emph{Traité des fonctions elliptiques et de leurs applications} (Treatise on elliptic functions and their applications) \cite{Halphen} carries the subtitle: \emph{Applications à la mécanique, à la physique, à la géodésie, à la géométrie et au calcul intégral} (Applications to mechanics, physics, geodesy,\index{geodesy} geometry and integral calculus).

Riemann adopted a physical  approach to functions of a complex variable. This point of view was new, compared to that of Cauchy, although both men reached simultaneously the characterization of conformal mappings in terms of the partial differential equations which are called the Cauchy--Riemann equations.\footnote{Riemann gives this characterization at the beginning of his doctoral dissertation\index{doctoral dissertation!Riemann}\index{Riemann!doctoral dissertation} \cite{Riemann-Grundlagen}, defended in 1851, and Cauchy in his papers \cite{Cauchy-Sur-1851} and \cite{Cauchy-1851-mono}, published the same year.} Riemann, just after establishing these equations, notes that the real and imaginary parts of such a function satisfy the Laplace equation.\index{Laplace equation}\index{equation!Laplace}  Ahlfors writes, in \cite{Ahlfors-dev} p. 4: ``Riemann virtually puts equality signs between two-dimensional potential theory and complex function theory."  We shall say more about Riemann's use of the Dirichlet principle,\index{Dirichlet principle}\index{principle!Dirichlet} in particular in his paper on Abelian functions \cite{Riemann-Abelian}, where he solves the question of the determination of a function of a complex variable by given conditions on the boundary and the discontinuity points. Klein, in his article on \emph{Riemann and his significance for the development of modern mathematics} (1895) \cite{Klein1923}, recalls the importance of potential theory and the influence of Dirichlet on Riemann. He writes:
\begin{quote}\small
It should also be observed that the theory of the potential, which in our day, owing to its importance in the theory of electricity and in other branches of physics, is quite universally known and used as an indispensable instrument of research, was at that time in its infancy. It is true that Green had written his fundamental memoir as early as 1828; but this paper remained for a long time almost unnoticed. In 1839, Gauss followed with his researches. As far as Germany is concerned, it is mainly due to the lectures of Dirichlet that the theory was farther developed and became known more generally; and this is where Riemann finds his base of operations.
\end{quote}
Among the works in potential theory that had a great impact later on, that of George Green,\footnote{George Green (1793--1841)\index{Green, George (1793--1841)} was a British mathematician and physicist who was completely self-taught. 
His father, also called George, was a baker, and the young George began working to earn his living at the age of five. He went to school for only one year, between the ages of 8 and 9. While he was working full-time in his father's mill, Green used the small amount of time that was left to him to study mathematics without the help of anybody else. On his own, Green became one of the main founders of potential theory. The word ``potential" was coined by him, although the notion existed before, e.g. in the works of Laplace\index{Laplace@de Laplace Pierre-Simon (1838--1749)}  and Poisson\index{Poisson, Sim\'eon Denis (1781--1840)} on hydrostatics. Besides, Green  developed mathematical theories of magnetism and  electricity that later on inspired the works of Maxwell\index{Maxwell, James Clerk (1831--1879)} and William Thomson  (later known as Lord Kelvin).\index{Thomson, William (Lord Kelvin) (1824--1907)} Green's father died one year after the publication of the his son's \emph{Essay on the Application of Mathematical Analysis to the Theories of Electricity and Magnetism}, printed at the author's own expense . In the meantime, George Sr. Green had gained some wealth, and what he left was sufficient for his son to put an end to his activities in the mill and to dedicate himself to mathematics.  At the recommendation of some influential acquaintance, Green was admitted at the University of Cambridge, as an undergraduate, in 1833, at the age of forty. The difficulties he encountered were not in catching up in the sciences, but in Greek and Latin. Green sat for the bachelor examination five years later, the same year as Sylvester, who was 21 years younger than him. 
During his relatively short career, Green wrote, besides the paper we mentioned above, several others, on optics, hydrodynamics, gravitation, and the theory of sound. He spent the last part of his life in Cambridge in isolation, addicted to alcohol. His work was rediscovered after his death by Thomson\index{Thomson, William (Lord Kelvin) (1824--1907)} and his ideas blossomed in physics and mathematics.}  mentioned by Klein in the last quote, is worth singling out because his author did it in isolation and never obtained, during his lifetime, the credit he deserves. Green, in 1828, gave the famous \emph{Green formula},\index{Green!formula} in his paper entitled \emph{An essay on the application of mathematical analysis to the theory of electricity and magnetism} \cite{Green}. This article contains the basis of what we call now Green's functions and Green's potential.\index{Green!potential} The introductory part of the essay emphasizes the role of a potential function, and this notion is then used in the setting of electricity and magnetism. The work also contains an early form of Green's Theorem\index{Green!theorem}\index{theorem!Green} (p. 11--12) which connects a line integral along a simple closed curve and the surface integral over the region bounded by that curve.\footnote{Basically, Green proved Stokes' theorem for surfaces embedded in 3-space.} There is an analogous theorem which relates volume and surface integrals contained in Riemann's 1851 inaugural dissertation.\index{doctoral dissertation!Riemann}\index{Riemann!doctoral dissertation} It might be noted that the result is stated (without proof) in an 1846 paper by Cauchy\index{Cauchy, Augustin-Louis (1789--1857)} \cite{Cauchy1846}. Cauchy\index{Cauchy, Augustin-Louis (1789--1857)} was also a physicist. We owe  him important works on hydrodynamics, elasticity, celestial mechanics and several other topics. But unlike Riemann, Cauchy's mathematical papers do not contain references to physics. Cauchy made a clear distinction between the methods of the two subjects.
In the introduction to the first volume of his famous \emph{Cours d'analyse de l'\'Ecole Polytechnique} \cite{Cauchy-cours}, we can read: 
\begin{quote}\small
Without any doubt, in the sciences which we call natural, the only method which is worth using with success consists in observing the facts and submitting later on the observations to calculus. But it would be a big mistake to think that we can only find certainty in geometric proofs or in the evidence of senses. [...] Let us cultivate with hard work the mathematical sciences, without intending to extend them beyond their domain; and let us not imagine that one can address history with formulae, neither giving as sanctions for morals, theorems from algebra or integral calculus.\footnote{Sans doute, dans les sciences qu'on nomme naturelles, la seule m\'ethode qu'on puisse employer avec succ\`es consiste à observer les faits et à soumettre ensuite les observations au calcul. Mais ce serait une erreur grave de penser qu'on ne trouve la certitude que dans les d\'emonstrations g\'eom\'etriques, ou dans le t\'emoignage des sens [...] Cultivons avec ardeur les sciences math\'ematiques, sans vouloir les \'etendre au-delà de leur domaine ; et n'allons pas nous imaginer qu'on puisse attaquer l'histoire avec des formules, ni donner pour sanction à la morale des th\'eor\`emes d'alg\`ebre ou de calcul int\'egral.}
\end{quote}

 Turning back to the Ancients, we quote a related phrase from Boook I of the  \emph{Nicomachean Ethics}; cf. \cite{Aristote-Nico} 1094b:\footnote{I thank M. Karbe for this reference.} 
 \begin{quote}\small [...] for it is the mark of an educated mind to expect that amount of exactness in each kind which the nature of the particular subject admits. It is equally unreasonable to accept merely probable conclusions from a mathematician and to demand strict demonstration from an orator. 
 \end{quote}

At the beginning of his dissertation, Riemann introduces the definition of conformality in terms of the existence of a complex derivative. From this point of view, a function $w$ of a complex variable $z$ is conformal if the derivative $\frac{dw}{dz}$ exists and is independent of the direction. This is equivalent to the infinitesimal notion of angle-preservation. As a matter of fact, conformality of maps in the sense of angle-preservation was already rooted in physics before Riemann. It is important to remember that the question of representing conformally the surface of a sphere onto the plane was already addressed by Hipparchus (second century B.C.), Ptolemy\index{Ptolemy, Claudius (c. 100--c. 170 A.D.)}  (first century A.D.), and certainly other Greek geometers and astronomers  in their work on spherical geometry and cartography, see \cite{Ptolemy-conformal} p. 405ff. We refer the reader to the recent surveys \cite{Pap-Geo}, \cite{Papa-Tissot}  regarding the relation between geography and conformal  and quasiconformal (in the sense of close-to-conformal) mappings. One may mention in particular Euler\index{Euler, Leonhard (1707--1783)} who studied general conformal maps from the sphere to the plane in his memoirs \cite{E490}, \cite{E491}, \cite{E492} which he wrote in relation with his work as a cartographer.\footnote{At the Academy of Sciences of Saint Petersburg, Euler, among his various duties, had the official charge of  cartographer and participated in the huge project of drawing maps of the Russian Empire.} In these memoirs, Euler expressed the conformality of projection maps from the sphere onto a Euclidean plane in terms of partial differential equations. Lambert in his paper \cite{Lambert1772}, also formulated problems concerning the projection of subsets of the sphere onto the Euclidean plane in terms of partial differential equations. Likewise, 
Lagrange\index{Lagrange, Joseph-Louis (1736--1813)}  used the notion of conformal map in his papers on cartography \cite{Lagrange1779}, and the same notion is inherent in  Gauss's work.\index{Gauss, Carl Friedrich (1777--1855)} The terminology ``isothermal coordinates" which the latter introduced, referring to a locally conformal map between a subdomain of the plane and a subdomain of the surface, indicates the relation with  physics. Riemann, in his dissertation, refers to an 1822 paper by Gauss, published in the \emph{Astronomische Abhandlungen} in 1825 (but written several years before).\footnote{The paper won a prize for a question proposed by the Copenhagen Royal Society of Sciences in 1822. The subject of the competition was:
 ``To represent the parts of a given surface onto another surface in such a way that the representation is similar to the original in its infinitesimal parts." A letter from Gauss to Schumacher dated July 5, 1816 shows that the solution was already known to Gauss at that time; cf. Gauss's \emph{Werke} Vol. 8, p. 371.} The title of that paper is \emph{Allgemeine Aufl\"osung der Aufgabe: die Theile einer gegebnen Fl\"ache auf einer andern gegebnen Fl\"ache so abzubilden, da\ss \  die Abbildung dem Abgebildeten in den kleinsten Theilen \"ahnlich wird} (General solution of the problem: to represent the parts of a given surface on another so that  the smallest parts of the representation shall be similar to the corresponding parts of the surface represented), a paper presented to a prize question proposed by the Royal Society of Sciences at Copenhagen, \cite{Gauss-Copenhagen}. In this paper, Gauss shows that every sufficiently small neighborhood  of a point in an arbitrary real-analytic surface can be mapped conformally onto a subset of the plane.\footnote{Gauss did not solve the problem of mapping conformally an arbitrary finite portion of the surface; this was one of the questions considered by Riemann. An English translation of Gauss's paper is published in \cite{Smith}, Volume 3.}

After recalling the definition of a conformal map, Riemann passes to the equivalent condition expressed in terms of partial differential equations. Here, we are given
  a function $f$ of a complex variable which is composed of two functions $u$ and $v$  of two real variables $x$ and $y$: 
\[f(x+iy)= u+iv.\]
 The functions $u$ and $v$ are differentiable and satisfy the Cauchy--Riemann equations.
They appear as \emph{potentials} in the space of the two variables $x$ and $y$. Klein, in his article on \emph{Riemann and his significance for the development of modern mathematics} writes (\cite{Klein1923} p. 168):
\begin{quote}\small
Riemann's method can be briefly characterized by saying that \emph{he applies to these parts $u$ and $v$ the principles of the theory of the potential}. In other words, his starting point lies in the domain of mathematical physics.
\end{quote}
In the same article (p. 170), after explaining some of Riemann's tools, Klein adds:  
\begin{quote}\small
All these new tools and methods, created by Riemann for the purpose of pure mathematics out of the physical intuition, have again proved of the greatest value for mathematical physics. Thus, for instance, we now always make use of Riemann's methods in treating the \emph{stationary} flow of a fluid within a two-dimensional region. A whole series of most interesting problems, formerly regarded as insolvable, had thus been solved completely. One of the best known problems of this kind is Helmholtz's\index{Helmholtz@von Helmholtz, Hermann Ludwig Ferdinand (1821--1894)} determination of the shape of a free liquid jet.
\end{quote}

Klein,\index{Klein, Felix (1849--1925)} who spent a significant part of his time advertising and explaining Riemann's ideas, completely adhered to his physical point of view. In 1882, he wrote a booklet entitled \emph{\"Uber Riemanns Theorie der algebraischen Funktionen und ihrer Integrale} (On Riemann's theory of algebraic functions and their integrals: A supplement to the usual treatises) \cite{Klein-Riemann} in which he explains the main ideas in Riemann's 1857 article on Abelian functions.
 This booklet is a redaction of part of a course that Klein gave in 1881 at  the University of Leipzig, and it had a great influence in making Riemann's ideas known.\footnote{Constance Reid reports on p. 178 of her biography of Hilbert \cite{Reid-Hilbert} that at a meeting of the G\"ottingen Scientific Society dedicated to the memory of Klein, held a few months after his death, Courant declared: ``If today we are able to build on the work of Riemann, it is thanks to Klein."} 
The excerpts we present here and later in this chapter are from the English translation \cite{Klein-Riemann}, published a few years after the German original. 

According to Klein,\index{Klein, Felix (1849--1925)}  the point of view on analytic functions based on the Cauchy--Riemann equations\index{Cauchy--Riemann equations}\index{equations!Cauchy--Riemann} is supported by physics. The first paragraph of his exposition   \cite{Klein-Riemann} (p. 1) is entitled \emph{Steady streaming in the plane as an interpretation of the functions of $x+iy$}.  He writes there: ``The physical interpretation of those functions of $x+iy$ which are dealt with in the following pages is  well known." He refers  to Maxwell's\index{Maxwell, James Clerk (1831--1879)} \emph{Treatise of electricity and magnetism} (1873), and he adds: ``So far as the intuitive treatment of the subject is concerned, his point of view is exactly that adopted in the text." Maxwell,\index{Maxwell, James Clerk (1831--1879)} at the beginning of the 1860s, developed a theory of electricity and magnetism and established the partial differential equations that carry his name, which describe the generation of electric and magnetic fields and the relation between them. In some sense, Maxwell's equations\index{Maxwell equations}\index{equations!Maxwell} are a generalization of the Cauchy--Riemann equations.\index{Cauchy--Riemann equations}\index{equations!Cauchy--Riemann} Klein\index{Klein, Felix (1849--1925)}  is among the first mathematicians who stressed this point.  After he states the  Cauchy--Riemann equations, Klein continues: 
\begin{quote}\small In these equations we take $u$ to be the \emph{velocity-potential}, so that $\frac{\partial v}{\partial y}, \frac{\partial u}{\partial x}$ are  the components of the velocity of a fluid moving parallel to the $xy$ plane [...] For the purposes of this interpretation it is of course indifferent of what nature we may imagine the fluid to be, but for many reasons it will be convenient to identify it here with the \emph{electric fluid}: $u$ is then proportional to the electrostatic potential which gives rise to the streaming, and the apparatus of experimental physics provide sufficient means for the production of many interesting systems of streaming.
\end{quote} Later in the text, in dealing with residues, Klein\index{Klein, Felix (1849--1925)}  writes: 
\begin{quote}\small The reason that the residue of $z_0$ must be equal and opposite to that of $z_1$ is now at once evident: the streaming is to be steady, hence the amount of electricity flowing at one point must be equal to that flowing out at the other.
\end{quote}
 In another report on Riemann's work, \cite{Klein1923} p. 175, Klein states: 
 \begin{quote}\small Riemann's treatment of the theory of function of complex variables, founded on the partial differential equation of  the potential, was intended by him to serve merely as an \emph{example} of the analogous treatment of all other physical problems that lead to partial differential equations, or to differential equations in general. [...] The execution of this programme which has since been considerably advanced in various directions, and which has in recent years been taken up with particular success by French geometers, amounts to nothing short of a \emph{systematic reconstruction of the methods of integration required in mechanics and in mathematical physics.}
 \end{quote}

 It is also interesting to recall that according to Klein \cite{Klein-development}, Riemann started studying Abelian functions because of their use in his research on galvanic currents.

An important element in Riemann's theory of functions of a complex variable is the so-called Dirichlet problem.\index{Dirichlet problem}\index{problem!Dirichlet}\footnote{Concerning the terminology, Klein writes in his \emph{Development of mathematics in the 19th century} ({Klein-development} p. 242 of the English translation): ``This is the first boundary value problem, which the French, unhistorical as they are, call the `Dirichlet problem': to determine a function $u$ if its boundary values and definite physically possible discontinuities are given -- there will be one and only one solution.}
Stated with a minimal amount of hypotheses, the problem, from the mathematical point of view, asks for the following: Given an open subset $\Omega$ of $\mathbb{R}^n$  and a continuous function $f$ defined on the boundary $\partial \Omega$ of $\Omega$, to find an extension of $f$ to $\Omega$ which is harmonic and continuous on the union $\Omega\cup\partial\Omega$. 
  The problem has more than one facet and there are several ways of dealing with it. Physicists consider that the problem has obviously a positive solution under very mild conditions, and that this solution is unique. In this setting, one thinks of the function $f$ on $\partial \Omega$ as a time-independent potential (electric, gravitational, etc.). Letting the system evolve, it will attain an equilibrium state, and the solution will necessarily satisfy a mean value property, that is, it will be harmonic. The harmonicity property is also formulated in terms of realizing the minimum of the energy functional 
  \begin{equation}\label{energy1}
\int\int\left( (\frac{\partial u}{\partial x})^2+(\frac{\partial u}{\partial y})^2\right)dxdy.
\end{equation}
 All these ideas were known to eighteenth century physicists and in fact, most of them can be traced back to Newton.

The ``Dirichlet principle"\index{Dirichlet principle}\index{principle!Dirichlet} is a method for solving the Dirichlet problem.\index{Dirichlet problem}\index{problem!Dirichlet} It is Riemann who coined the term. The principle is based on an assertion he took for granted, namely, that an infimum of the energy functional is attained.  This infimum is necessarily harmonic and froms a harmonic function $u$. 
    Riemann used the Dirichlet principle\index{Dirichlet principle}\index{principle!Dirichlet} to construct analytic functions, not only on the disc, but on an arbitrary Riemann surface, after cutting it along a system of arcs so that it becomes simply-connected. Riemann also used the Dirichlet principle\index{Dirichlet principle}\index{principle!Dirichlet} at other places in his doctoral dissertation\index{doctoral dissertation!Riemann}\index{Riemann!doctoral dissertation} and in his paper on Abelian functions.\footnote{The name ``Dirichlet principle"\index{Dirichlet principle}\index{principle!Dirichlet} is used in the paper \cite{Riemann-Abelian} on Abelian functions (\S III, IV, Preliminaries), but not in the doctoral dissertation \cite{Riemann-Grundlagen}. Riemann, in his existence proof of meromorphic functions on general Riemann surfaces, defined these functions by their real parts, which are harmonic functions, using this principle. He also used in his proof of the Riemann Mapping Theorem. In fact, it is well known today that the Riemann Mapping Theorem,\index{Riemann mapping theorem}\index{theorem!Riemann mapping} the existence of meromorphic functions, and  the Dirichlet problem,\index{Dirichlet!problem} are all equivalent.   
Riemann writes in \S III of the Preliminary section of his paper on Abelian functions \cite{Riemann-Abelian} that in the study of integrals of algebraic functions and their inverses, one can use a principle  which Dirichlet used several years before in his lectures on the forces that act by the inverse of the square of the distance, to solve of a problem related to a function of three variables satisfying the Laplace equation.\index{Laplace equation}\index{equation!Laplace} He adds that Dirichlet was probably inspired by an analogous idea of Gauss. In fact, Gauss used such a principle in his  1839 paper \cite{Gauss-inverse}. He assumed there without proof that for a given constant potential distribution, an equilibrium state is  attained and is unique and corresponds to the minimum of the energy. It is possible that Riemann chose to call this principle the ``Dirichlet principle"\index{Dirichlet principle}\index{principle!Dirichlet} out of faithfulness to the mathematician from whom he learned most. Klein writes in \cite{Klein-development}: ``Riemann was bound to Dirichlet by the strong inner sympathy of a like mode of thought. Dirichlet loved to make things clear to himself in an intuitive substrate; along with this he would give acute, logical analyses of foundational questions and would avoid long computations as much as possible. His manner suited Riemann, who adopted it and worked according to Dirichlet's methods."} 
 At the time Riemann appealed to the Dirichlet problem,\index{Dirichlet problem}\index{problem!Dirichlet} several other eminent mathematicians used an analogous principle, in physics and in mathematics.  This includes Laplace, Fourier and Poisson in France, Green, Thomson\index{Thomson, William (Lord Kelvin) (1824--1907)} and Stokes in England,  and Gauss in Germany. Helmholtz\index{Helmholtz@von Helmholtz, Hermann Ludwig Ferdinand (1821--1894)} used this principle in his work on acoustics \cite{Helm-Crelle}.  Riemann's use of the Dirichlet principle\index{Dirichlet principle}\index{principle!Dirichlet} was criticized by Weierstrass\index{Weierstrass, Karl (1815--1897)} \cite{W14-1870}. Klein writes in his \emph{Development of mathematics in the 19th century} ({Klein-development} p. 248 of the English translation): 
 \begin{quote}\small
 
  With this attack by Weierstrass on Dirichlet's principle,\index{Dirichlet principle}\index{Dirichlet principle}\index{principle!Dirichlet} the evidence to which Dirichlet, and after him, Riemann, had appealed, became fragile [...]  The majority of mathematicians turned away from Riemann; they had no confidence in the existence theorems, which Weierstrass's critique had robbed of their mathematical supports.
  
 The physicists took yet another position: they rejected Weierstrass's critique. Helmholtz,\index{Helmholtz@von Helmholtz, Hermann Ludwig Ferdinand (1821--1894)} whom I once asked about this, told me: ``For us physicists the Dirichlet principle\index{Dirichlet principle}\index{principle!Dirichlet} remains a proof." Thus he evidently distinguished between proofs for mathematicians and physicists; in any case, it is a general fact that physicists are little troubled by the fine points of mathematics -- for them the ``evidence" is sufficient.
 \end{quote}

    The mathematicians' doubts concerning Riemann's use of the Dirichlet principle\index{Dirichlet principle}\index{principle!Dirichlet} were removed only several years later. We refer the reader to \cite{Monna} for the details of this interesting story. 
%

 In the preface to his booklet \cite{Klein-Riemann}, (p.  IX), Klein\index{Klein, Felix (1849--1925)}  writes:  
\begin{quote}\small 
[...] there are certain physical considerations which have been lately developed, although restricted to simpler cases, from various points of view.\footnote{[Klein's footnote] Cf. C. Neumann, Math. Ann. t. X, pp. 569-771. Kirchoff, Berl. Monatsber., 1875, pp. 487-497. T\"opler, Pogg. Ann. t.  CLX., pp. 375-388.} I have not hesitated to take these physical conceptions as the starting point of my presentation. Riemann, as we know, used Dirichlet's Principle in their place in his writings. But I have no doubt that he started from precisely those physical problems, and then, in order to give what was physically evident the support of mathematical reasoning, he afterwards substituted Dirichlet's Principle. 
\end{quote}
Klein adds:
\begin{quote}\small
Anyone who clearly understands the conditions under which Riemann worked in G\"ottingen, anyone who has followed Riemann's speculations as they have come down to us, partly in fragments, will, I think, share my opinion. However that may be, the physical method seemed the true one for my purpose. For it is well known that Dirichlet's principle is not sufficient for the actual foundation of the theorems to be established; moreover, the heuristic element, which to me was all-important, is brought out far more prominently by the physical method. Hence the constant introduction of intuitive considerations, where a proof by analysis would not have been difficult and might have been simpler, hence also the repeated illustration of general results by examples and figures.
\end{quote}
Poincar\'e\index{Poincar\'e, Henri (1854--1912)}  was, like Riemann, a pre-eminent representative of a philosophical tradition of thought in geometry and physics which was invoked at the outset, a tradition combining mathematical, physical, and philosophical thinking. In his booklet  { \emph{La valeur de la science} (The value of science)  \cite{Poincare1911} (1905), commenting on Klein's method, Poincar\'e writes: 
   \begin{quote}\small
[...] On the contrary, look at Mr. Klein: He is studying one of the most abstract questions in the theory of functions; namely, to know whether on a given Riemann surface there always exists a function admitting given singularities: for instance, two logarithmic singular points with equal residues of opposite signs. What does the famous German geometer do? He replaces his Riemann surface by a metal surface whose electric conductivity varies according to certain rules. He puts the two logarithmic points in contact with the two poles of a battery. The electric current must necessarily pass, and the way this current is distributed on the surface defines a function whose singularities are the ones prescribed by the statement.\footnote{[...] Voyez au contraire M. Klein: il étudie une des questions les plus abstraites de la théorie des fonctions ; il s'agit de savoir si sur une surface de Riemann donnée, il existe toujours une fonction admettant des singularités données : par exemple, deux points singuliers logarithmiques avec des résidus égaux et de signe contraire. Que fait le célèbre géomètre allemand ? Il remplace sa surface de Riemann par une surface métallique dont la conductibilité électrique varie suivant certaines lois. Il met les deux points logarithmiques en communication avec les deux p\^oles d'une pile. Il faudra bien que le courant passe, et la fa\c con dont ce courant sera distribué sur la surface définira une fonction dont les singularités seront précisément celles qui sont prévues par l'énoncé.}
\end{quote}
  
To end this section, let us mention some of the numerous applications of Riemann surfaces in modern mathematical physics.  

One of the major applications of the theory of Riemann surfaces in physics is the Atiyah-Singer index theorem.\index{Atiyah-Singer theorem}\index{theorem!Atiyah-Singer}  This theorem, obtained in 1963, gives an information on the dimension of the space of solutions of a differential operator (the \emph{analytical index}) in terms of topology (the \emph{topological degree}). The theorem is used in the theory of the Einstein equation, the instanton equation, the Dirac operator, etc. It is considered as a vast generalization of the classical version of the theorem of Riemann--Roch\index{Riemann--Roch theorem}\index{theorem!Riemann--Roch}, which is an equality, half of which contained in Riemann's paper on Abelian functions \cite{Riemann-Abelian},  and the other half in the dissertation of his student Roch\index{Roch, Gustav (1839-1866)} \cite{Roch}. (See \cite{AAF} in this volume for a review of this theorem.)

We also mention string theory,\index{string theory} in which (0-dimensional) particles of physics are replaced by (1-dimensional) strings. This theory was developed as a framework that would hopefully solve some problems that cannot be handled by the theory of relativity. At some point (and it still is, albeit with a more skepticism on physical and mathematical grounds) string theory was considered as a possible theory for the unification of the fundamental forces in nature: gravitation and quantum theory, including electromagnetism -- another attempt to realize Riemann's long-life insight. In this theory, one follows the history of a closed string, that is, a closed loop in 3-space. While it propagates, such a loop sweeps out a surface. For reasons that have to do with the consistency of the theory, the surface turns out to be equipped with a 1-dimensional complex structure, that is, it is a Riemann surface. If the string does not interact with anything else, then the swept-out surface is a cylinder, but in general, the string, under some interaction, splits into two other strings, which join again, etc. creating a Riemann surface of higher connectivity (Figures \ref{string1} and  \ref{string2}). Seen from very large distances, strings look like ordinary particles, they have mass and charge, but they can also vibrate. This vibration leads to a hypothetical quantum mechanical particle called graviton, which is supposedly responsible for the gravitational force. It is in this sense that string theory is a theory of quantum gravity.

%
%
%

 Riemann surfaces are at the basis of conformal field theories (CFT),\index{conformal field theory} in which one associates to a marked Riemann surfaces a vector space satisfying certain natural axioms. These surfaces also appear as a major ingredient in the topological quantum field theories (TQFT) developed by Witten and others, which are based on sets of axioms that provide  functors from a certain category of cobordisms to the category of vector spaces (Segal and Atiyah gave such sets of axioms). TQFTs lead to results in physics (relativity, quantum gravity, etc.) and at the same time to results in mathematics, where they provide quantum invariants of 3-manifolds. They have applications in symplectic geometry, representation theory of Lie groups and algebraic geometry, in particular in the study of moduli spaces of holomorphic vector bundles over Riemann surfaces. One may also mention that the famous geometric Langlands correspondence is based on the theory of Riemann surfaces. Stated  loosely, in the geometric Langlands correspondence  one assigns to each rank $n$ holomorphic vector bundle  with a holomorphic connection on a complex algebraic curve, a Hecke eigensheaf on the moduli space of rank $n$ holomorphic vector bundles on that curve, cf. \cite{Frenkel}.

 Riemann surfaces are also the main ingredients in the theory of 
Higgs bundles. These objects arose in the study made by Nigel Hitchin of the self-duality equation on a Riemann surface. From the physical point of view, Higgs bundles describe  particles like the Higgs boson. Conversely, the   physical methods of Higgs bundle theory are used in the study of  moduli spaces of  representation of surface groups. Hitchin's motivation arose from his work  done in the 1970s  with Atiyah, Drinfield and Manin on the so-called instanton equation,\index{instanton equation}\index{equation!instanton} another theory combining in an essential way mathematics and physics \cite{AH}.
 
Finally, let us mention that Riemann surfaces are used in biology, cf. the recent survey \cite{Penner}.

\section{Riemann's memoir on trigonometric series\index{Riemann! habilitation text}\index{habilitation text!Riemann}\index{trigonometric series}}\label{s:trigo}

The habilitation degree, which was required in Germany in order to hold a university teaching position, involved two presentations: the \emph{Habilitationsschrift},\index{Riemann! habilitation text}\index{habilitation text!Riemann} a written original work on a specialized subject, and the \emph{Habilitationsvortrag},\index{Riemann! habilitation lecture}\index{habilitation lecture!Riemann} a lecture on a subject chosen by the university council. 
The present section is devoted to 
 Riemann's  \emph{Habilitationsschrift}\index{Habilitationsschrift}  \cite{Riemann-Trigo}. We shall discuss his Habilitationsvortrag\index{Riemann! habilitation lecture}\index{habilitation lecture!Riemann} in the next one.

Riemann's  Habilitationsschrift is entitled \emph{\"Uber die Darstellbarkeit einer Function durch eine trigonometrische Reihe} (On the representability of a function by a trigonometric series).\index{trigonometric series}\index{Riemann! habilitation text}\index{habilitation text!Riemann} It is generally considered that Riemann worked on it during thirty months. He presented it to the university in December 1853.   
 About this work, in a letter to Hou\"el, dated  March 18, 1873, and quoted in \cite{Dugac1973}, Darboux writes: ``This memoir of Riemann is a masterpiece which is similar to these old paintings of which some small parts in full light make you regret what time has destroyed or what the author has neglected."\footnote{Ce m\'emoire de Riemann est un chef-d'\oe uvre semblable à ces vieux tableaux dont quelques parties en pleine lumi\`ere vous font regretter ce que le temps a d\'etruit ou ce que l'auteur a n\'eglig\'e.}

 This theory of trigonometric series finds its origin\index{trigonometric series} in eighteenth century physics, more precisely, in the introduction by d'Alembert,\index{Alembert@d'Alembert, Jean le Rond (1717--1783)} in 1747, of the vibrating string equation (also called the wave equation). To understand the context of Riemann's contribution, it might be useful to recall a few key events in the history of the subject. This theory expanded very slowly, and it was eventually put on firm bases in the nineteenth century, mainly by Joseph Fourier,\index{Fourier, Joseph (1768--1830)}  while he was working on another problem arising from physics, namely, heat diffusion. In the meantime, many pre-eminent mathematicians and physicists worked on trigonometric series,\index{trigonometric series} and we shall mention a few of them. Furthermore, the work done during the first decades after the introduction to the vibrating string equation gave rise to one of the most passionate controversies in the history of mathematics and physics whose scope was larger than the subject of trigonometric series, and we shall say a few words about it. The controversy involved Euler\index{Euler, Leonhard (1707--1783)}, d'Alembert,\index{Alembert@d'Alembert, Jean le Rond (1717--1783)}  
 Lagrange,\index{Lagrange, Joseph-Louis (1736--1813)}  Daniel Bernoulli\index{Bernoulli, Daniel (1700--1782)} and other major scientists. In particular, a  quarrel between Euler and d'Alembert lasted from 1748 until 1783 (the year both of them died).  Later on, a dispute concerning the same subject broke out between Fourier and Poisson.\index{Poisson, Sim\'eon Denis (1781--1840)} The question was about the ``continuity" of the functions representing the solutions. This dispute is thoroughly discussed in the introductory part of volume IV of Series A of Euler's \emph{Opera omnia} \cite{O. IV A 5}, a volume containing the correspondence between Euler\index{Euler, Leonhard (1707--1783)} and d'Alembert.\index{Alembert@d'Alembert, Jean le Rond (1717--1783)} A comprehensive  survey of this controversy is also made in \cite{Jehel} and in Chapter 1 of the present volume \cite{Papa-Euler}.

 When d'Alembert discovered the vibrating string equation,  Euler\index{Euler, Leonhard (1707--1783)} immediately became interested. He had already been dealing with partial differential equations for several years. In fact, he started working, around the year 1735, on partial differential equations and their applications in geometry and physics.  Furthermore, the theory of sound was one of his favorite subjects.\footnote{Euler writes in his memoir \cite{LE_1} that ``the most sublime research that scientists successfully undertook these days  is in all respects without question that of propagation of sound." [La plus sublime recherche que les géomètres aient entreprise de nos jours avec succès est sans contredit à tous égards celle de la propagation du son.] We also recall that the subject of Euler's\index{Euler, Leonhard (1707--1783)} first published memoir is the theory of sound \cite{Euler_dissertatio-J}.}\index{sound!theory of}  This subject was not new, and, in fact,  it is worth recalling that the physics of vibrating strings was one of the main problems studied by the Pythagoreans, back in the sixth century B.C. Indeed, most of the ancient biographers of Pythagoras\index{Pythagoras (c. 580 B.C.-- c. 495 B.C.)} describe his experiments on pitch production, cf. \cite{Iamblichus}. For a recent scholarship on Pythagoras and the early Pythagoreans, the reader may consult \cite{Zhumd}.

The heart of the controversy on the vibrating string lies in the question of the clarification of the notion of function, more precisely, the nature of the functions that are solutions of the partial differential equation representing the vibration of a string.\index{string vibration}\index{vibrating string}

 We discuss this matter in Chapter 1 of the present 
volume \cite{Papa-Euler}. Instead, we make here an excursion to the origin of the theory of sound production,\index{sound!theory of}  in order to make fully clear that Riemann's investigations on trigonometric functions and integration theory originate in physics.

The first part of Riemann's \emph{Habilitationsschrift}\index{Riemann! habilitation text}\index{habilitation text!Riemann} is a historical report on the representation of a function by a trigonometric series, and in fact, it is motivated by the theory of theory of the vibrating string. In a letter to his father, written in the autumn of 1852, Riemann says that he learned the historical details from Dirichlet, who explained them to him in a two-hour session. (The letter is reproduced in Riemann's Collected Works, \cite{Riemann-Gesammelte} p. 578.) 
Riemann starts his historical survey by recalling that this subject is important for physics:
\begin{quote}\small
Trigonometric series, which are given this name by Fourier, that is, series of the form
\[a_1\sin x+a_2\sin 2x+a_3\sin 3x+\ldots\]
\[+\frac{1}{2}b_0+b_1\cos x+b_2\cos 2x+b_3\cos 3x+\ldots\]
play a substantial role in the part of mathematics where we encounter functions which are completely arbitrary. We also have reasons to say that the progress of this part of mathematics, which is so important for physics, has been subject to a more precise knowledge of the nature of these series.
\end{quote}
Riemann's excursion in history is divided into three periods, and we shall say a few words about each period.

The first period is concerned with the controversies that arose concerning the notion of function which led to the question of representing arbitrary functions by a trigonometric series.\index{trigonometric series} D'Alembert,\index{Alembert@d'Alembert, Jean le Rond (1717--1783)} in  1747, wrote two papers, which were published in the Memoirs of the Berlin Academy and under the titles \emph{Recherches sur la courbe que forme une corde tendue mise en vibration} (Researches on the curve that is formed by a stretched vibrating string) \cite{Al_recherches}  and \emph{Suite des recherches sur la courbe que forme une corde tendue mise en vibration} (Sequel to the researches on the curve that is formed by a stretched vibrating string) \cite{Al_suite} .
 In these memoirs, d'Alembert, relying on the fundamental principle of dynamics, gave the partial differential equation that represents the motion of a point on a vibrating string\index{vibrating string!equation}\index{equation!vibrating string} subject to small vibrations:\index{string vibration}\index{vibrating string} 
  \begin{equation} \label{eq:wave}
  \frac{\partial^2y}{\partial t^2}=\alpha^2\frac{\partial^2y}{\partial x^2}.
  \end{equation}
  Here, $\alpha$ is a constant and $y$ is the oscillation of the string, a function of time, $t$, and distance along the string, $x$. 
  In the same memoir, d'Alembert wrote the  first general solution to the problem, with the given boundary conditions,  in the form
\[y(x,t)=\frac{1}{2}\left(\phi(x+\alpha t)+\phi(x-\alpha t)\right),\] where $\phi$ is an ``arbitrary" periodic function whose period is the double of  the length of the string.   D'Alembert\index{Alembert@d'Alembert, Jean le Rond (1717--1783)} used a method he attributes to Euler\index{Euler, Leonhard (1707--1783)} for the integration of partial  differential equations.  The problem was to give a meaning to the adjective ``arbitrary," and this is where the more basic question of \emph{What is a function?} was raised.
 
It is natural to assume that the solution of d'Alembert's vibration equation should be (twice) differentiable, since the equation involves second partial derivatives, and this is what d'Alembert did. Euler was not of the same opinion. The reason he gave is physical, namely, that one can give a non-smooth initial form to the string (for example a curve with corners) which is being pinched, therefore the function that represents the shape of the string could be quite arbitrary. This implies that the solution may be arbitrary. Euler published his remarks in his  memoir \emph{Sur la vibration des cordes} (On the vibration of strings)\footnote{Euler wrote two versions, one in Latin and one in French, the French version bearing the mention ``Translated from the Latin."} \cite{Euler_vib-corde} in which he reviews d'Alembert's work on the wave equation. These remarks introduced some doubts concerning the work of d'Alembert, who wrote a new memoir on the same subject, in which he confirms his ideas, \emph{Addition aux recherches sur la courbe que forme une corde tendue mise en vibration} (Addition to the researches on the curve formed by a stretched vibrating string)
\cite{Al_addition}.

Another pre-eminent scientist who became involved in these questions was Daniel Bernoulli, who was primarily a physicist. Before talking about his contribution to the subject, one should recall that Brook Taylor,\index{Taylor, Brook (1685-1731)} 
 in his memoir \emph{De motu nervi 
tensi} (On the motion of a tense string) \cite{Taylor1713}  (1713)
and later in his treatise  \emph{Methodus incrementorum directa et inversa} (Direct and Indirect Methods of Incrementation) \cite{Taylor1715}; first edition 1715, noted that a trigonometric function like $f(x)=\sin x$ represents a periodic phenomenon, a wave.  In his work on the subject, Taylor was motivated by music theory. 
Bernoulli came out with  a formula of the  form
\begin{equation}\label{Ber} \sum_{n=1}^{\infty}\sin \frac{n\pi x}{l}\cos \frac{n\pi \alpha t}{l}.
\end{equation} Like Taylor, he was motivated by music. In fact, among all the scientists of the Bernoulli family, Daniel was the most inclined  towards physics. His intuition concerning Formula (\ref{Ber}) originates in the fact known to all music theorists that the string vibration\index{string vibration}\index{vibrating string} produces, together with the fundamental pitch, an infinite sequence of harmonics. It is interesting to read some excerpts of Bernoulli's writings on this subject.  In his memoir \emph{ Mémoire   sur   les   vibrations   des   cordes   d'une   épaisseur   inégale} (Memoir on the vibrations of a string of uneven width), he writes \cite{Ber1767} (p. 173):
 \begin{quote}\small
 I showed furthermore in the Berlin Memoirs that the vibrations of various orders, however one takes them, may coexist in one and the same string, without disturbing each other in any way, these various kinds of coexisting vibration being absolutely independent of each other. Hence this multiplicity of harmonic sounds which we hear at the same time and with one and the same string. If all the modes of vibration start at the same instant, it may happen that the first vibration of first order, the second vibration of second order, the third vibration of third order, etc. terminate at the same moment. 
This is in some sense an apparent synchronism which is nothing less than general, since there are infinitely many vibrations which do not terminate at the same instant.\footnote{J'ai démontré de plus dans les Mémoires de Berlin, que les vibrations de différents ordres,
quels   qu'on   les   prenne,   peuvent   coexister   dans   une   seule   et   même   corde,   sans   se   troubler   en aucune façon, ces différentes espèces de vibration coexistantes étant absolument indépendantes les unes des autres. 
De là cette pluralité de sons harmoniques qu'on entend à la fois d'une seule  et 
même   corde.   Si   toutes   espèces   de   vibration   commencent   au   même   instant,   il   arrivera   que   la
premi\`ere vibration du premier ordre, la seconde vibration du second ordre, la troisième vibration
du   troisième   ordre   etc.   finiront   au   même   instant.   C'est   là   un   synchronisme   apparent   dans   un
certain sens, et qui n'est rien moins que général, puisqu'il y a une infinité d'autres vibrations qui
ne finissent pas au même instant.}

\end{quote}

In his \emph{Réflexions et éclaircissements sur les   nouvelles   vibrations des cordes} (Reflections and clarifications on the new vibrations of strings), Bernoulli writes  (\cite{Bernoulli1755a} p. 152--153):

\begin{quote}\small

Indeed, all musicians agree that a long pinched string gives at the same time, besides its fundamental tone, other tones which are much more acute; most of all they will notice the mixture of the twelfth and the minor sixteenth: in case they don't notice as much distinctly the octave and the double octave, it is only because of the very big resemblance of these two tones with the fundamental. This is an evident proof that there could occur in one and the same string a mixture of several sorts of Taylorian vibrations at the same time. In the same manner, we hear in the sound of large bells a mixture of different tones. If we hold by the middle a steel stick, and if we hit it, we hear at the same time a confused mixture of several tones, which, when appreciated by a skilled musician, turn out to be extremely inharmonious, in such a way that a combination of vibrations is formed, which never start and finish at the same moment, except by a happenstance: hence we see that the harmony of sounds, which we hear  at the same time in one sonorous body, is not essential to that material, and should not serve as a principle for systems in music. Air is not free of this multiplicity of coexisting sounds: it often happens that one extracts two different sounds from a pipe; but the best proof of how much the various air waves may prevent each other is that we hear distinctly every part of a concert, and that all the waves due to these different parts are formed from the same mass of air without disturbing each other, very much like light rays entering in a dark room from a small hole do not disturb each other.\footnote{Effectivement tous les Musiciens conviennent, qu'une longue corde pincée donne en même temps,
outre son ton fondamental, d'autres tons beaucoup plus aigus ; ils remarquent surtout le m\'elange
de   la   douzième   et   de   la   dix-septième   majeure   :   s'ils   ne   remarquent   pas   aussi   distinctement l'octave et la double octave, ce n'est qu'à cause de la trop grande ressemblance de ces deux
tons  avec  le  ton fondamental.  Voilà une  preuve  évidente,  qu'il peut  se   faire  dans  une  seule  et 
même   corde   un   mélange   de   plusieurs   sortes   de   vibrations   Tayloriennes   à   la   fois.   On   entend
pareillement dans le son des grosses cloches un mélange de tons différents. Si l'on tient par le
milieu une verge d'acier, et qu'on la frappe, on entend à la fois un mélange confus de plusieurs
tons, lesquels étant appréciés par un habile Musicien se trouvent extr\^emement désharmonieux, de
sorte qu'il se forme un concours de vibrations, qui ne commencent et ne finissent jamais dans un
même instant, sinon par un grand hazard : d'où l'on voit que l'harmonie des sons, qu'on entend
dans une même corps sonore à la fois, n'est pas essentielle à cette matière, et ne doit pas servir de
principe   pour   les   systèmes   de   Musique.   L'air   n'est   pas   exempt   de   cette   multiplicité   de   sons coexistants : il arrive souvent qu'on tire deux sons différents à la fois d'un tuyau ; mais, ce qui prouve   le   mieux,   combien   peu   les   différentes   ondulations   de   l'air   s'entre-empêchent,   est   qu'on entend distinctement toutes les parties d'un concert, et que toutes les ondulations causées par ces différentes   parties   se   forment   dans   la   même   masse   d'air   sans   se   troubler   mutuellement,   tout comme   les   rayons   de   la   lumière,   qui   entrent   dans   une   chambre   obscure   à   travers   une   petite ouverture, ne se troublent point.}
\end{quote}
 
 In the same memoir, Bernoulli writes (\cite{Bernoulli1755a} p. 151):
\begin{quote}\small

My conclusion is that every sonorous body contains essentially an infinity of sounds, and infinitely corresponding ways of performing their regular vibrations; finally, that in each different way of vibrating the variations in the parts of the sonorous body are formed in a different way.\footnote{Ma conclusion est, que tous les corps sonores renferment en puissance une infinité de sons, et une infinité de manières correspondantes de faire leurs vibrations régulières ; enfin, que dans
chaque différentes espèce de vibrations les inflexions des parties du corps sonore se font d'une
manière différente.}
\end{quote}

 In the same memoir, Bernoulli writes 
 (\cite{Bernoulli1755a} p. 148):
\begin{quote}\small
[...] without any less esteem for the calculations of Messrs d'Alembert and Euler, which certainly include the most profound and exquisite things that analysis contains; but which show at the same time that an abstract analysis, which we follow without any synthetic examination of the proposed question, may be surprising rather than enlightening for us. It seems to me that one had only to be attentive to the nature of simple vibrations of a string in order to foresee without any calculation everything these two geometers found by the most tricky and abstract calculations with which an analytical mind has been instructed.\footnote{[...] je n'en estime pas
moins les calculs de Mrs. d'Alembert et Euler, qui renferment certainement tout ce que l'Analyse
peut avoir de plus profond et de plus sublime ; mais qui montrent en même temps, qu'une analyse
abstraite, qu'on écoute sans aucun examen synthétique de la question proposée, est sujette à nous
surprendre plutôt qu'à nous éclairer. Il me semble à moi, qu'il n'y avait qu'à faire attention à la
nature des vibrations simples des cordes, pour prévoir sans aucun calcul tout ce que ces deux
grands géomètres ont trouvé par les calculs les plus épineux et les plus abstraits, dont l'esprit
analytique se soit encore avisé.}

\end{quote}

 It may seem surprising that Euler,\index{Euler, Leonhard (1707--1783)} who was as much involved in music theory than Bernoulli -- he had even corresponded with Rameau on overtones back in 1752 (see \cite{Euler-Hermann}, Vol. II) --  did not state this idea before, especially that Euler had already heavily manipulated infinite series.  It is also a fact that the techniques of trigonometric series are quite different from those of power series.

 When Daniel Bernoulli suggested that an ``arbitrary" function defined on a finite interval can be expanded as a convergent trigonometric series, several basic questions appeared at the forefront of research:  
 \begin{enumerate}
 \item What is the meaning of such an infinite sum, that is, in what sense does it converge?
 \item
In what sense functions possess    trigonometric series expansions, and how can such a result be proved.
 \item What is an ``arbitrary" function (a question that had been thoroughly investigated without reaching any definite conclusion), and more precisely, is there a definition of an arbitrary function such that it coincides with functions expressible by such an infinite sum?  
 \end{enumerate}
 
The second period of 
   Riemann's historical report is dominated by Joseph Fourier (1768--1830)\index{Fourier, Joseph (1768--1830)} who, in his \emph{Th\'eorie analytique de la chaleur} (Analytic theory of heat) \cite{Fourier} (1822), developed the theory of trigonometric series\index{trigonometric series}, while he was studying the heat equation. Let us say a few words about Fourier's treatise.
   
   The introduction (Discours préliminaire) of this treatise is interesting. It concerns the importance of heat in our universe. Fourier\index{Fourier, Joseph (1768--1830)} writes, at the beginning of that introduction (p. i): 
   \begin{quote}\small Heat, like gravity, penetrates all substances of the universe, its rays occupies all parts of space. The goal of our work is to present the mathematical laws that govern this element. From now on, this theory will constitute one of the most important branches of general physics.\footnote{La chaleur pénètre, comme la gravité, toutes les substances de l'univers, ses rayons occupent toutes les parties de l'espace. Le but de notre ouvrage est d'exposer les lois mathématiques que suit cet élément. Cette théorie formera désormais une des branches les plus importantes de la physique générale.}
   \end{quote} Fourier\index{Fourier, Joseph (1768--1830)} then mentions the works of Archimedes, Galileo and Newton,\index{Newton, Isaac (1643--1727)} and he comments on the importance of the effect of sun rays on every element of the living world. It is interesting to note that Archimedes, Galileo and Newton are again mentioned, together, in the introduction to Riemann's habilitation lecture\index{Riemann! habilitation lecture}\index{habilitation lecture!Riemann} 1854, which we consider in \S \ref{s:vortrag}. 
   
    After that, he arrives at the mathematical principles of that theory. The problem, which turned out to be very difficult to solve, is stated very clearly (\cite{Fourier} p.  2):
      \begin{quote}\small  When heat is unevenly distributed between the various points of a solid mass, it tends to an equilibrium position, and it slowly passes from the overheated parts to the ones which are less heated. At the same time, it dissipates through the surface, and it gets lost in the ambient space or in void. This tendency towards a uniform distribution, and this spontaneous emission which takes place at the surface of bodies, causes a continuous change in the temperature at the various points. The question of the propagation of heat consists in determining what is the temperature at each point of a body at a given time, assuming that the initial temperatures are known.\footnote{Lorsque la chaleur est inégalement distribuée entre les différents points d'une masse solide, elle tend à se mettre en équilibre, et passe lentement des parties les plus échauffées dans celles qui le sont moins ; en m\^eme temps elle se dissipe par la surface, et se perd dans le milieu ou dans le vide. Cette tendance à une distribution uniforme, et cette émission spontanée qui s'opère à la surface des corps, changent continuellement la température des différents points. La question de la propagation de la chaleur consiste à déterminer quelle est la température de chaque point d'un corps à un instant donné, en supposant que les températures initiales sont connues.}
          \end{quote}

 In Chapter II of his treatise, Fourier\index{Fourier, Joseph (1768--1830)} establishes (p. 136) the  equation which is known nowadays as the ``heat equation":
\[\frac{\partial u}{\partial t}=k^2\left( \frac{\partial^2u}{\partial x^2}+\frac{\partial^2u}{\partial y^2}+\frac{\partial^2u}{\partial z^2}\right).\]  
In Chapter III, he gives the solution of this equation in the form of a trigonometric series.\index{trigonometric series} On page 243, he writes: ``The preceding analysis gave us the way to develop an arbitrary function as a series of sines and cosines of multiple arcs." Then he announces that he will apply these  results to some particular cases which show up in physics, as solutions of partial differential equations. On page 249, he considers the problem of the vibrating string, and he declares that the principles he established solve the difficulties that are inherent in the analysis done by Daniel Bernoulli.\index{Bernoulli, Daniel (1700--1782)} He recalls that the latter gave a solution that assumes that an arbitrary function can be developed as a trigonometric series,\index{trigonometric series} but that of all the proofs that were proposed of this fact, the most complete is the one where we can determine the coefficients of such a function. 
This is precisely what Fourier does. For a given trigonometric series,\index{trigonometric series}\index{Riemann! habilitation text}\index{habilitation text!Riemann}
\[f(x)=a_1\sin x+a_2\sin 2x+\ldots -\frac{1}{2}b_0+b_1\cos x+b_2\cos 2x+\ldots,
\] Fourier\index{Fourier, Joseph (1768--1830)} provides the (now well-known)  integral formula for  the coefficients:
\[a_n=\frac{1}{\pi}\int_{-\pi}^{\pi} f(x)\sin nx dx\]
and 
\[b_n=\frac{1}{\pi}\int_{-\pi}^{\pi} f(x)\cos nx dx.\]

Picard, in the series of three lectures on the history of analysis that he gave in America \cite{Picard-developpement}, says  (p. 7) that these integral formulae were known to Euler, who mentioned them incidentally. In the same lectures, Picard insists on the fact that Fourier had an audacious method which involved the solution of an infinite number of first-order equations, with an infinite number of unknowns.

Fourier also showed that his theory can be applied to functions which may have discontinuities\footnote{The word ``discontinuity" is understood here in the modern sense of the word, and not in the sense of Euler. Cf. the explanation in Chapter 1 of the present volume \cite{Papa-Euler}.}.  It is also good to recall that Fourier\index{Fourier, Joseph (1768--1830)} stated, back in 1807 (the paper was published in 1808, \cite{Fourier1808}), the fact that a function, given graphically in an arbitrary manner, may be expressed by a trigonometric series.  

Riemann reports in his memoir on trigonometric functions that at that time Lagrange\index{Lagrange, Joseph-Louis (1736--1813)}  vigorously rejected Fourier's assertion. He also recalls the rivalry between Fourier\index{Fourier, Joseph (1768--1830)} and Poisson,\index{Poisson, Sim\'eon Denis (1781--1840)} and that the latter took the defense of Lagrange. Riemann analyzes some passages from Lagrange's work, and he concludes by repeating that Fourier was the first to understand exactly and completely the nature of trigonometric series. He adds that after Fourier's work, these series appear in several ways in mathematical physics, as representations of arbitrary functions. He declares that in each particular case one was able to prove that the Fourier series  indeed converges to the value of  the function, but that it took a long time before such an important theorem was proved in full generality. He recalls that in 1826 Cauchy attempted a proof of that result using complex numbers, in a memoir of the Academy of Sciences (t. VI, p. 603), but that this proof is incomplete, as was shown by Dirichlet.\footnote{Riemann refers to \cite{Dirichlet-Crelle}.}  Riemann declares that he completed Cauchy's proof in his inaugural dissertation.\index{doctoral dissertation!Riemann}\index{Riemann!doctoral dissertation}
 
 The third section of the historical part of Riemann's \emph{Habilitationsschrift}\index{Riemann! habilitation text}\index{habilitation text!Riemann} concerns the work of Dirichlet.\index{Dirichlet (Lejeune), Peter Gustav (1805--1859)} The latter, whom we already mentioned several times and who had been one of Riemann's teachers in Berlin, gave a necessary and sufficient condition under which a periodic function can be expanded as a trigonometric series (\cite{Dirichlet-Crelle}, 1829). In the same memoir, Dirichlet\index{Dirichlet (Lejeune), Peter Gustav (1805--1859)} obtained the general theorem concerning the convergence of Fourier series after he pointed out some errors in Cauchy's proof of that result.   Riemann, in his \emph{Habilitationsschrift},\index{Riemann! habilitation text}\index{habilitation text!Riemann} considers that it is Dirichlet\index{Dirichlet (Lejeune), Peter Gustav (1805--1859)} who closed the controversy. 
 He declares that the latter, in a publication which appeared in 1829 in Crelle's journal (t. IV),\footnote{This is Dirichlet's article \cite{Dirichlet1829}.} gave a ``very rigorous" theory of representation by trigonometric series\index{trigonometric series} of general periodic functions under the hypothesis that they are integrable, that they do not have infinitely many maxima or minima, and that at the points of discontinuity, the value of the function is the arithmetic mean of its left and right limits. Dirichlet\index{Dirichlet (Lejeune), Peter Gustav (1805--1859)} left open the converse: given a function that does not satisfy the first two conditions (the third one must obviously be satisfied), under what conditions can it be represented by a trigonometric series?\index{trigonometric series} 
      This is one of the questions that Riemann solved in his habilitation memoir. For more details, the interested reader is referred to the exposition in Chapter 1 of the present volume \cite{Papa-Euler}.
 
   In conclusion, the question of the meaning of a function started with physics: the vibration equation discovered by d'Alembert,  and it ended again with physics: the study of heat, by Fourier, and the discovery of Fourier series, which are extensively used in mathematical physics. Picard writes in his historical survey \cite{Picard-developpement} that the development of a function as a series is a remarkable example of the intimate solidarity that unites at certain points pure analysis and applied mathematics.

   \section{Riemann's Habilitationsvortrag\index{Riemann! habilitation lecture}\index{habilitation lecture!Riemann} 1854 --   Space and Matter} \label{s:vortrag} 
   
Riemann's public lecture, his Habilitationsvortrag,\index{Riemann! habilitation lecture}\index{habilitation lecture!Riemann}  \emph{\"Uber die Hypothesen, welche der Geometrie zu Grunde liegen},  which  was the final requirement before he was allowed to teach at the university level, was delivered on June 10, 1854.  This lecture marks the birth of modern differential geometry. It is a difficult text, involving --  like for other writings of Riemann -- mathematics, physics and philosophy.  It was commented on by many philosophers and scientists, and translated several times into other languages.\footnote{There are English translations by M. Spivak in his \emph{Comprehensive Introduction to Differential Geometry} (\cite{Spivak},  Volume 2, pp. 132--153) and by H. S. White in Smith's \emph{Source book mathematics} (\cite{Smith}, p. 411--425), and probably others. J. Jost's edition \cite{Riemann-Jost} contains Cayley's translation. Italians translations were made by E. Betti, and E. Betti, and a French one by J. Ho\"uel.} 
 The earliest translation is probably the one that Clifford\index{Clifford, William Kingdon (1845--1879)} made in 1873 for the journal \emph{Nature} \cite{Clifford}. This translation is generally considered as too literal. It is nevertheless interesting because Clifford was at the same time a mathematician, a physicist and a philosopher. He was knowledgeable in the philosophical issues raised by Riemann and he was familiar with the specialized language of philosophy that the latter used. 
  
In a broad sense, the subject of the investigation is geometry, space and the relation between them. The discussion takes place at several levels, starting from the foundations of geometry: Riemann mentions the axioms at the beginning of his essay. He introduces several kinds of spaces and the notion of ``manifoldness"\index{manifoldness} which we shall discuss below. He mentions in particular discrete and continuity manifoldnesses, infinitesimal and large-scale properties properties, the ambient physical space,\index{space} mathematical $n$-dimensional spaces and $n$-uply extended magnitudes.
He declares that the propositions of geometry cannot be derived  from the general notion of magnitude (the word is taken in the Aristotelean sense), and that the properties which distinguish (physical) space from other conceivable triply extended magnitudes are only to be deduced from experience. 

Several authors commented on Riemann's dissertation, and we shall make a few remarks on them below. 
Riemann's lecture has three parts (We use a slight modification of Clifford's original headlines in \emph{Nature}):
\begin{enumerate}
\item The notion of an $n$-tuply extended magnitude.
\item Measure-relations of which a manifoldness of $n$ dimensions is susceptible, on the assumption that lines have a length independent of position, and consequently that every line may be measured by every other.
\item Applications to space.
\end{enumerate}
Roughly speaking, the first part is philosophical, the second one is mathematical, and the third one deals with applications to physics. But to some extent philosophy and physics  are present in the three parts. A detailed explanation of all these notions would take us too far, and it is also known that several points in this essay are very cryptic. Heinrich Weber, Hermann Weyl and many other pre-eminent mathematicians tried to uncover their meaning. Weyl, who had a great devotion for Riemann, edited the Habilitationsvortrag\index{Riemann! habilitation lecture}\index{habilitation lecture!Riemann} in 1919, \cite{Weyl1919}, together with a commentary, making the link with relativity theory. One of the main features of the local geometry conceived by Riemann is that it is well suited to the study of gravity and more general fields in physics. Relativity theory, which encompasses the largest part of modern physics, relies in a crucial way on the notions introduced by Riemann.

  From the purely mathematical point of view, the most important contribution of the Habilitationsvortrag\index{Riemann! habilitation lecture}\index{habilitation lecture!Riemann} is that it sets the bases of what we call today Riemannian geometry, with the introduction of the curvature tensor and its consequences, including several results such as the fact that the homogeneity (with the inherent notion of transformation group) corresponds to constant curvature. This new geometry can be considered as a far-reaching generalization of Gauss's work on the intrinsic geometry of surfaces, and at the same time it is a generalization of the non-Euclidean geometry (of constant curvature) discovered by Lobachevsky, Bolyai and Gauss, in the few decades that preceded Riemann. Furthermore, Riemann sets in this memoir the bases of several developments made in several directions by Clifford, Christoffel, Bianchi, Ricci, Beltrami, Levi-Civita, \'Elie Cartan,  Einstein and many others.

 It is known that the full importance of the Habilitationsvortrag\index{Riemann! habilitation lecture}\index{habilitation lecture!Riemann} was not recognized in the first years after it was delivered.
A report by Dedekind on the mathematical content of the memoir was published only in 1868.\footnote{\emph{Abhandlungen der K\"oniglichen Gesellschaft der Wissenschaften zu G\"ottingen, 13.}}  But it is also known that Gauss, who, as Riemann's mentor, was present at the lecture, expressed his complete  satisfaction with it. This reaction to Riemann's Habilitationsvortrag\index{Riemann! habilitation lecture}\index{habilitation lecture!Riemann} is described in Dedekind's biography of Riemann published in the Collected Works \cite{Riemann-Gesammelte}.  Gauss's praise was certainly a cause for Riemann's own contentedness, because Gauss was known to be sparing with compliments.

In the introduction to his dissertation,\index{doctoral dissertation!Riemann}\index{Riemann!doctoral dissertation} Riemann declares that he is unexperienced in philosophical questions, and that in preparing the lecture he could rely only on some remarks that Gauss made in his second paper on biquadratic residues and in his Jubilee-book, and some philosophical researches of Herbart.  Riemann discusses the question of space and that of manifoldness and its specialization. This specialization may be either continuous  or discrete. In the chapter \cite{Plot} contained in the present volume, Plotnitsly emphasizes that Riemann speaks of discrete manifolds, and then says that, rather than space itself, it is ``the reality underlying space" that may be discrete.\footnote{It may be useful to note that in modern physics, spacetime\index{Minkowski spacetime}\index{spacetime} is studied in its both characters, discrete (e.g., in lattice gauge theories, which are often considered as mathematical discrete approximations) or as continuous (e.g., in general relativity).} ``Manifolds"\index{manifold} as we intend them today are particular cases of manifoldnesses. Riemann writes:
\begin{quote}\small
Manifoldnesses in which, as in the plane and in space, the line-element may be reduced to the form $\sqrt{ \sum dx^2 }$, are therefore only a particular case of the manifoldnesses to be here investigated; they require a special name, and therefore these manifoldnesses in which the square of the line-element may be expressed as the sum of the squares of complete differentials I will call \emph{flat}. 

\end{quote}

Riemann declares that notions whose specializations form a continuous manifoldness are the positions of perceived objects (die Orte der Sinnengegenst\"ande) and colors. It is conceivable that Riemann, in his mention of colors, refers to the fact that one can continuously move from a color to another one, a color being characterized by the proportions of red, green and violet it contains.  This makes color a three-dimensional manifoldness. Weyl writes, in his \emph{Space, time, matter}, that Riemann's reference to color\index{color} ``is confirmed by Maxwell's\index{Maxwell, James Clerk (1831--1879)} familiar construction of the color triangle" (\cite{Weyl-STM}, p. 84 of the English translation). There are also writings of Helmholt\index{Helmholtz@von Helmholtz, Hermann Ludwig Ferdinand (1821--1894)}z and Thomas Young\index{Young, Thomas (1773--1829)} on this matter.
 There are particular portions of a manifoldness called \emph{quanta},\footnote{In Aristotle's\index{Aristotle (384--322 B.C.)} language, the latin word \emph{quantum} is used as the translation of the word $\pi o \sigma \acute{o} \nu$ (quantity, or ``quantified thing"), which is one of Aristotle's categories.\index{Categories (Aristotle)} In his \emph{Metaphysics}, Aristotle\index{Aristotle (384--322 B.C.)} mentions four types of change: of  substance, \emph{quale},  \emph{quantum}, or place. (Metaphysics, 1069b9-13).} whose nature is different from that which is characterized by the discrete and the continuous. ``Their comparison with regard to quantity is accomplished in the case of discrete magnitudes by counting, in the case of continuous by measuring." The notions of measurement and  of dimension are discussed in the Aristotelean style. The relation between measurement and the axioms of geometry is said to be fundamental, and in some sense this question concerns the relation between the axioms of geometry and the reality of space, that is, between mathematical and empirical truth. Riemann writes:
 \begin{quote}\small
 Either therefore the reality which underlies space must form
a discrete manifoldness, or we must seek the ground of its metric relations outside it, in binding forces which act upon it.\footnote{Es muss also entweder das dem Raume zu Grunde liegende Wirkliche eine discrete Mannigfaltigkeit bilden,  oder der Grund der Massverhältnisse ausserhalb, in darauf wirkenden bindenden Kräften, gesucht werden.}
 \end{quote}

  It is important to recall Riemann's words. He declares that geometry depends at the same time on axioms and on observational and experimental physics.  He considers that classical geometry, with the first principles and axioms that it assumes and the connections between them, does not lead anywhere, because ``we do not perceive the necessity of these connections." What is missing is a notion of ``multiply-extended magnitude" (mehrfach ausgedehnte Gr\"osse), a notion which makes space a particular ``triply-extended magnitude." He proposes that the properties that distinguish space from other conceivable  triply-extended magnitudes be deduced from experience. 
In particular, the space that Riemann talks about, although built from undefined notions and axioms connecting them, is not the space of traditional geometry. He suggests that this space should reflect the material world around us. He formulates the problem of finding the ``simplest matters of fact from which the metric relations (Massverh\"altnisse) of space may be determined." He declares that these matters of fact are ``not of necessity, but only of empirical certainty."  They are the ``hypotheses" that are referred to in the title of the dissertation. He says that he will investigate the ``probability" of these matters of fact, ``within the limits of observation," and see whether they may be extended ``beyond the limits of observation, both on the side of the infinitely great and of the infinitely small."   The most important among these matters of fact is related to the work of Euclid. The geometry that Riemann will construct will be Euclidean at the infinitesimal level.\footnote{It may be worth recalling that Lobachevsky, in his various works on non-Euclidean geometry that he started in the late 1820s, systematically checked that the formulae that obtained in his new geometry give, at the infinitesimal level, the Euclidean formulae. See e.g. \cite{Lobachevsky-Pan} p. 31.}
We already noted by the way that the notion of ``infinitely small" is treated in several works of the Ancient Greeks. The same notion is thoroughly discussed by Galileo Galilei\index{Galilei, Galileo (1564--1642)} in his \emph{Discorsi; First day}, to whom Riemann refers in his habilitation, though in a different context.

The rest of the Habilitationsvortrag\index{Riemann! habilitation lecture}\index{habilitation lecture!Riemann} is a development of the ideas expressed in the introduction. The first part  concerns the notion of $n$-dimensional Mannigfaltigkeit.\index{Mannigfaltigkeit (manifoldness)}\index{manifoldness}  This term is sometimes translated into English by ``manifoldness."\index{Mannigfaltigkeit (manifoldness)}\index{manifoldness}
  Riemann also talks about a ``multiply extended magnitude." This is an ancestor to the mathematical notion of manifold. But the meaning of Mannigfaltigkeit, in Riemann's terminology, and that of the mathematical notion of manifold, as it is used today, do not coincide, even though in German the word Mannigfaltigkeit is used for ``manifold."\footnote{Jost, in \cite{Riemann-Jost},  tries to sort out this complex terminology. He writes on p. 29: ``The English of Clifford may appear somewhat old-fashioned for a modern reader. For instance, he writes `manifoldness' instead of the simpler modern translation `manifold' of Riemann's term `Mannigfaltigkeit.' But Riemann's German sounds likewise  somewhat old-fashioned, and for that matter, `manifoldness' is the more accurate translation of Riemann's term. In any case, for historical reasons, I have selected that translation here."} There are discrete\index{manifoldness!discrete}\index{discrete manifoldness} and continuous manifoldnesses,\index{manifoldness!continuous}\index{continuous manifoldness} and there are manifoldnesses\index{manifoldness} which are not mathematical. Riemann says that notions with
specializations to discrete manifoldness are very common, but that, by contrast, there are very few  notions whose specialization form a continuous manifoldness.\index{continuous manifoldness}\index{manifoldness!continuous} We note incidentally that Poincar\'e pondered this terminology. In a letter to the mathematician and historian of mathematics Gustav Enestr\"om,\index{Enestr\"om, Gustav (1852--1923)}  dated November 19, 1883 (cf. \cite{Poincare-corresp} p. 143), he writes: 
 \begin{quote}\small I prefer the translation of \emph{Mannigfaltigkeit}  by multiplicity, because the two words have the same etymological meaning. The word set is more adapted to the \emph{Mannigfaltigkeiten}  considered by M. Cantor and which are discrete. It would be less adapted to those which I consider and which are discontinuous. What is the opinion of M. Mittag on this matter?\footnote{Je pr\'ef\`ere la traduction de \emph{Mannigfaltigkeit} par multiplicit\'e, car les deux mots ont même sens \'etymologique. Le mot ensemble convient bien aux \emph{Mannigfaltigkeiten} envisag\'es par M. Cantor et qui sont discr\`etes ; il conviendrait moins à celles que je consid\`ere et qui sont discontinues. Qu'en pense M. Mittag\index{Mittag-Leffler, G\"osta (1846--1927)} à ce sujet ?} 
 \end{quote} 
 Enestr\"om responds, on November 23: 
 \begin{quote}\small
 Mr. Mittag-Leffler thinks that you may be right, and, consequently, one should prefer the word multiplicity.\footnote{M. Mittag-Leffler pense que vous pouvez avoir raison, et que, par cons\'equent, il faut pr\'ef\'erer le mot multiplicit\'e.}
 \end{quote} In fact, Poincar\'e  used the word multiplicity in its French form (``multiplicit\'e") to denote Riemann's moduli space. 

In the treatise \emph{Th\'eorie des fonctions alg\'ebriques de deux variables ind\'ependantes} (Theory of algebraic functions of two independent variables) by Picard\index{Picard, Charles \'Emile (1856--1941)} and Simart\index{Simart, Georges (1846--1921)} \cite{Picard-Simart} on which we report in Chapter 8 of the present volume \cite{Papa-Riemann3}, the authors use interchangeably the words ``vari\'et\'e," ``multiplicit\'e" and ``continuum" to denote ``a certain \emph{continuous set} of points depending on a number of parameters which is equal to the dimension of this variety or continuum." (p. 20) It is interesting to note that Grothendieck,\index{Grothendieck, Alexandre (1928--2014)} in his \emph{Esquisse de programme} (A sketch of a program) \cite{Gro-esquisse}, uses the same word. We refer the reader to the chapters \cite{Ohshika} by Ohshika and \cite{Plot} by Plotnitsky in the present volume for further discussion of the notion of Mannigfaltigkeit in relation with manifolds.  
%
%
%
%
The second part of the habilitation lecture,\index{Riemann! habilitation lecture}\index{habilitation lecture!Riemann} which concerns metric relations (Maassverh\"altnisse), is more mathematically-oriented. It contains, condensed in six or seven pages, the foundations of Riemannian geometry. 
The exposition contains a minimum amount of formulae. In fact, there are essentially two formulae. The first one gives the line element in (Euclidean) ``space" as a square root of squares of differentials of the coordinates:
\begin{equation}\label{Riemann:Pytha} 
ds^2=\sum dx_i^2
.\end{equation}
This formula is an ``infinitesimal Euclidean Pythagorean theorem." expresses the fact that at the infinitesimal level the metric is Euclidean.
The second formula gives the line element in a curved space:
\[\frac{1}{1+\frac{\alpha}{4}\sum x^2}\sqrt{\sum dx^2}
\]
where $\alpha$ denotes, in Riemann's notation, the curvature.
As is well known, this formula gives the one of the Poincar\'e metric\index{Poincar\'e metric}\index{metric!Poincar\'e} of the disc in the case where the curvature constant is negative.

 After Riemann gives the general expression of the infinitesimal line element as the square root of a quadratic form, the curvature representing a deviation from flatness, he states that to know curvature at a point in a manifoldness\index{manifoldness} of dimension $n$, it is sufficient to know it in $n(n-1)/2$ surface directions. 
He notes that if the length of a line element is independent from its position (that is, the group of motions acts transitively), then the space must have constant curvature.  The third part, called ``Applications to space," contains in particular Riemann's famous discussion of the difference between \emph{unbounded} and \emph{infinite extent}. We quote him again:\begin{quote}\small
In the extension of space-construction to the infinitely great, we must distinguish between \emph{unboundedness} and \emph{infinite extent}, the former belongs to the extent relations, the latter to the measure-relations. That space is an unbounded three-fold manifoldness,\index{manifoldness} is an assumption which is developed by every conception of the outer world; according to which every instant the region of real perception is completed and the possible positions of a sought object are constructed, and which by these applications is for ever confirming itself. The unboundedness of space possesses in this way a greater empirical certainty than any external experience. But its infinite extent by no means follows from this; on the other hand if we assume independence of bodies from position, and therefore ascribe to space constant curvature, it must necessarily be finite provided this curvature has ever so small a positive value. If we prolong all the geodesics starting in a given surface-element, we should obtain an unbounded surface of constant curvature, i.e., a surface which in a \emph{flat} manifoldness of three dimensions would take the form of a sphere, and consequently be finite.
\end{quote}
This is a famous passage for which Riemann's name is associated with the geometry of the sphere (constant positive curvature).\index{geometry!spherical}\index{spherical geometry} It has been commented on by mathematicians, and also by philosophers, and it is related to the question of whether the universe has a spherical shape or not.
Again, we can quote related texts from Greek antiquity, e.g. from Empedocles\index{Empedocles of Agrigentum (c. 492--c. 432 B.C.)} concerning the universe as a round ``boundless" sphere, of which only a few fragments remain \cite{Empedocles}: 
\begin{quote}\small
The Sphere on every side the boundless same,
\\
Exultant in surrounding solitude.
\end{quote}

One may also quote Plato,\index{Plato (5th--4th c. B.C.)} who considers,  for philosophical reasons,  in the \emph{Timaeus} (\cite{Plato}, 33b), that the universe is spherical, and hence, bounded:
\begin{quote}\small
  He wrought it into a round, in the shape of a sphere, equidistant in all directions from the center to the extremities, which of all shapes is the most perfect and the most self-similar, since he deemed that the similar is infinitely fairer than the dissimilar. And on the outside round about, it was all made smooth with great exactness, and that for many reasons. 
  \end{quote}
Other Greek philosophers considered that the universe is infinite. This is an endless discussion.

Riemann wanted his (Riemannian)\index{Riemannian geometry}\index{geometry!Riemannian} geometry to represent at the same time the large-scale and, most of all, the small-scale geometries of space. 
The progress made in mechanics in the preceding centuries, he says, is due to the invention of the infinitesimal calculus and to the simple principles discovered by Archimedes,\index{Archimedes (c. 287 B.C.--c. 212 B.C.)} Galileo\index{Galilei, Galileo (1564--1642)} and Newton.\index{Newton, Isaac (1643--1727)} The natural sciences, are still in want of \emph{simple} principles. Let us quote Riemann again:
 \begin{quote}\small
 The questions about the infinitely great are for the interpretation of nature useless questions. But this is not the case with the questions about the infinitely small. It is upon the exactness with which we follow phenomena into the infinitely small that our knowledge of their causal relations essentially depends. The progress of recent centuries in the knowledge of mechanics depends almost entirely on the exactness of the construction which has become possible through the invention of the infinitesimal calculus, and through the simple principles discovered by Archimedes, Galileo, and Newton, and used by modern physics. 
 \end{quote}
 
The synopsis of the Habilitationsvortrag\index{Riemann! habilitation lecture}\index{habilitation lecture!Riemann} (Clifford's translation \cite{Clifford}) ends with two questions.
\begin{enumerate}
\item How far is the validity of these empirical determinations probable beyond the limits of observations towards the infinitely great?
\item How far towards the infinitely small? Connection of this question with the interpretation of nature.
\end{enumerate}

It is interesting to put again in parallel Riemann's writings with some texts of Aristotle\index{Aristotle (384--322 B.C.)} on related matters, and there are many of them. We choose an excerpt from the beginning of the treatise \emph{On the Heavens} \cite{Aristotle-Heavens}. It concerns at the same time magnitude, continuum, divisibility, infinity, dimension, infinitesimals, and the importance of these questions for understanding nature:
\begin{quote}\small
The science which has to do with nature clearly concerns itself for the most part with bodies and magnitudes and their properties and movements, but also with the principles of this sort of substance, as many as they may be. For of things constituted by nature some are bodies and magnitudes, some possess body and magnitude, and some are principles of things which possess these. Now a continuum is that which is divisible into parts always capable of subdivision, and a body is that which is every way divisible. A magnitude if divisible one way is a line, if two ways a surface, and if three a body. [...] All magnitudes, then, which are divisible are also continuous. Whether we can also say that whatever is continuous is divisible does not yet, on our present grounds, appear. [...] 
The question as to the nature of the whole, whether it is infinite in size or limited in its total mass, is a matter for subsequent inquiry. [...]
 This being clear, we must go on to consider the questions which remain. First, is there an infinite body, as the majority of the ancient philosophers thought, or is this an impossibility? The decision of this question, either way, is not unimportant, but rather all-important, to our search for the truth. It is this problem which has practically always been the source of the differences of those who have written about nature as a whole.  So it has been and so it must be; since the least initial deviation from the truth is multiplied later a thousandfold. Admit, for instance, the existence of a minimum magnitude, and you will find that the minimum which you have introduced, small as it is, causes the greatest truths of mathematics to totter.
 \end{quote}
Concerning the particular notion of infinite, we choose two texts from Aristotle's  \emph{Physics} \cite{Aristotle-Physics}. We stress on the fact that even though the Greek philosophers, represented by Aristotle, did not formulate an axiomatic (in a  mathematical sense) notion of infinity as we do it today, one should not underestimate the importance of the fact that they considered this notion as a fundamental philosophical notion, they asked many questions around it and about its role, and they also regarded it as central in physics and in mathematics.

The first text we choose is from Book III of Aristotle's \emph{Physics}:
\begin{quote}\small
[...] But on the other hand to suppose that the infinite does not exist in any way leads
obviously to many impossible consequences: there will be a beginning and an end of
time, a magnitude will not be divisible into magnitudes, number will not be infinite.
If, then, in view of the above considerations, neither alternative seems possible, an
arbiter must be called in; and clearly there is a sense in which the infinite exists and
another in which it does not.
We must keep in mind that the word ``is" means either what potentially is or what fully
is. Further, a thing is infinite either by addition or by division.
Now, as we have seen, magnitude is not actually infinite. But by division it is infinite.
(There is no difficulty in refuting the theory of indivisible lines.) The alternative then
remains that the infinite has a potential existence.
\end{quote}

The second text is from Book V of the  \emph{Physics}:
\begin{quote}\small Now it is impossible that the infinite should be a thing which is in itself infinite, separable from sensible objects. If the infinite is neither a magnitude nor an aggregate, but is itself a substance and not an accident, it will be indivisible; for 
the divisible must be either a magnitude or an aggregate. But if indivisible, then not infinite, except in the way in which the voice is invisible. But this is not the way in which it is used by those who say that the infinite exists, nor that in which 
we are investigating it, namely as that which cannot be gone through. But if the infinite is accidental, it would not be, qua infinite, an element in things, any more than the invisible would be an element of speech, though the voice is invisible. 
 
Further, how can the infinite be itself something, unless both number and magnitude, of which it is an essential attribute, exist in that way? If they are not 
substances, a fortiori the infinite is not. 
 
It is plain, too, that the infinite cannot be an actual thing and a substance and 
principle. For any part of it that is taken will be infinite, if it has parts; for to be 
infinite and the infinite are the same, if it is a substance and not predicated of a 
subject. Hence it will be either indivisible or divisible into infinites.  But the same 
thing cannot be many infinites. (Yet just as part of air is air, so a part of the infinite 
would be infinite, if it is supposed to be a substance and principle.) Therefore the 
infinite must be without parts and indivisible. But this cannot be true of what is 
infinite in fulfillment; for it must be a definite quantity. 
 
 Belief in the existence of the infinite comes mainly from five considerations:
From the nature of time -- for it is infinite; From the division of magnitudes -- for
the mathematicians also use the infinite [...]
\end{quote}

Finally, let us note that in Book IV of the \emph{Physics} \cite{Aristotle-Physics} there is a long discussion about time, with its relation to measure and change: 
\begin{quote}\small
As to what time is or what is its nature, the traditional accounts give us as little light. [...] It is evident, then, that time is neither movement nor independent of movement. 
We must take this as our starting-point and try to discover -- since we wish to know what time is -- what exactly it has to do with movement. 
\end{quote}

  Closer to us, another mathematician-philosopher who was fascinated by infinity is Blaise Pascal. He wrote on this theme, in his mathematical and philosophical writings. From his \emph{Pens\'ees} \cite{Pascal-P}, we read: 
\begin{quote}\small
Unity added to infinity adds nothing to it, any more than does one foot
added to infinite length. The finite is annihilated in presence of the infinite, and
becomes pure nothingness."\footnote{L'unit\'e jointe à l'infini ne l'augmente de rien, non plus qu'un pied à une mesure infinie, le fini s'an\'eantit en pr\'esence de l'infini et devient un pur n\'eant.} 
\end{quote}
\begin{quote}\small
Our soul has been cast into the body, where it finds
number, time and dimension. It reasons thereupon, and calls it nature, necessity,
and can believe nothing else.\footnote{Notre \^ame est jet\'ee dans le corps où elle trouve nombre, temps, dimensions, elle raisonne l\`a-dessus et appelle cela nature, n\'ecessit\'e, et ne peut croire autre chose.}
\end{quote}

  \begin{quote}\small
The eternal silence of these infinite spaces terrifies me.\footnote{Le silence \'eternel de ces espaces infinis m'effraye.}
\end{quote}

One should also talk about modern physics, where the same kind of questions are still the basic ones: What is space? What is time?  What physical theories describe at the same time the macroscopic and the microscopic worlds? What are the relations between these worlds? How do we pass between the discrete and the continuous?  

Riemann's last sentence in the Habilitationsvortrag\index{Riemann! habilitation lecture}\index{habilitation lecture!Riemann} shows that he modestly considered that in his work, he did not make any significant advance in the direction of physics: 
\begin{quote}\small
This leads us into the domain of another science, of physics, into which the object of this work does not allow us to go today.
\end{quote}
 Higher-dimensional spaces, from the mathematical point of view were surely considered before Riemann. But for the first time, Riemann's major achievement was to introduce on these spaces a geometry that was necessary for the development of modern physics. The physical theories of superstrings and supergravity need ten or eleven dimensions.
 The spacetime of special relativity -- Minkowski's spacetime\index{Minkowski spacetime}\index{spacetime} --  is a four-dimensional manifold equipped with a structure that generalizes the one that Riemann considered. In an address to the 80th Assembly of German Natural Scientists and Physicians, (Sep 21, 1908), Minkowski\index{Minkowski, Hermann (1864--1909)} declares (cf. \cite{Minkowski1908}): 
   \begin{quote}\small
   The views of space and time which I wish to lay before you have sprung from the soil of experimental physics, and therein lies their strength. They are radical. Henceforth, space by itself, and time by itself, are doomed to fade away into mere shadows, and only a kind of union of the two will preserve an independent reality. 
   \end{quote}
In this setting, Riemann's formula (\ref{Riemann:Pytha}) is replaced by the formula
\begin{equation}\label{element-spacetime}
ds^2=c^2dt^2-dx^2-dy^2-dz^2
\end{equation}
where $t$ is the time component, $(x,y,z)$ the space components and $c$ the velocity of light. The geometry of Minkowski spacetime\index{Minkowski spacetime} is included in the setting of semi-Riemannian geometry,\index{semi-Riemannian geometry}\index{geometry!semi-Riemannian} a geometry in which the metric tensor is not necessarily positive-definite. This incorporates in the theory the fact that particles cannot move at a speed which is larger than the speed of light. But the basic features that Riemann conceived are there. In general relativity,\index{general relativity} the metric tensor is an expression of the gravitational potential, in the trend of  Riemann's ideas.

Let us now mention some comments by various authors (especially physicists) on the Habilitationsvortrag.\index{Riemann! habilitation lecture}\index{habilitation lecture!Riemann} We quote Clifford and Weyl whom we already mentioned.

On February 21, 1870, Clifford\index{Clifford, William Kingdon (1845--1879)} presented a paper to the Cambridge Philosophical Society whose title is \emph{On the space theory of matter} \cite{Clifford1870}, in which he stressed the relation of the new geometry with physics. It is interesting to read the abstract of that paper, for it gives quite a good idea of its physical background. Clifford\index{Clifford, William Kingdon (1845--1879)} writes:
\begin{quote}\small
Riemann has shown that as there are different kinds of lines and surfaces, so there are different kinds of spaces of three dimensions; and that we can only find out by experience to which of these kinds the space in which we live belongs. In particular, the axioms of plane geometry are true within the limits of experiment on the surface of a sheet of paper, and yet we know that the sheet is really covered with a number of small ridges and furrows, upon which (the total curvature not being zero) these axioms are not true. Similarly, he says that although the axioms of solid geometry are true within the limits of experiment for finite portions of our space, yet we have no reason to conclude that they are true for very small portions; and if any help can be got thereby for the explanation of physical phenomena, we may have reason to conclude that they are not true for very small portions of space.

I wish here to indicate a manner in which these speculations may be applied to the investigation of physical phenomena. I hold in fact

(1) That small portions of space \emph{are} in fact of a nature analogous to little hills on a surface which is on the average flat; namely, that the ordinary laws of geometry are not valid in them.

(2) That this property of being curved or distorted is continually being passed on from one portion of space to another after the manner of a wave.

(3) That this variation of the curvature of space is what really happens in that phenomenon which we call the \emph{motion of matter}, whether ponderable or etherial. 

(4) That in the physical world nothing else takes place but this variation, subject (possibly) to the law of continuity.

I am endeavoring in a general way to explain the laws of double refraction on this hypothesis, but have not yet arrived at any results sufficiently decisive to be communicated. 
\end{quote}
 
It is not superfluous to recall that Gauss,\index{Gauss, Carl Friedrich (1777--1855)} who was Riemann's mentor, was also interested in the philosophical implications of the new discoveries of geometry. In a letter dated March 6, 1832 (see \cite{Staeckel} and Gauss's Collected Works Vol. VI \cite{Gauss-collected}), Gauss writes to his friend Wolfgang Bolyai\index{Bolyai, Wolfgang (1775--1856)} that Kant\index{Kant, Immanuel (1724--1804)}  was wrong in declaring that space is \emph{only the form}\footnote{Gauss's emphasis.}  of our intuition. These remarks are made amid a discussion on non-Eucldiean geometry. In the same letter, Gauss refers to an article on the subject that he published in the G\"ottingische Gelehre Anzeigen, in 1831. This article is contained in Volume II of Gauss's \emph{Werke}. Gauss criticizes an argument, which is independent of non-Euclidean geometry, which Kant\index{Kant, Immanuel (1724--1804)} gave in support of his assumption (and his proof) that space is only a form of our exterior intuition. The argument is in Kant's \emph{Prolegomena zu einer jeden k\"unftigen Metaphysik, die als Wissenschaft wird auftreten k\"onnen} (Prolegomena to any future metaphysics that will be able to present itself as a science) \cite{Kant} \S 13, and it is based on the existence of symmetries. Gauss's position was that, on the contrary, space has a real significance, independent of our mode of intuition. An excerpt on space of Kant's inaugural dissertation -- in fact an excerpt concerned by Gauss's critic -- is quoted in Chapter 1 of the present volume (in the section concerning space).

 The questions of space and of time remained among the major preoccupations of Kant.\index{Kant, Immanuel (1724--1804)}  They are developed in particular in his habilitation \cite{Kant-Dissertation} and in his \emph{Critik der reinen Vernunft} (Critique of pure reason) \cite{Kant-Critik} (1781), which is one of the most influential philosophical works ever written. 
In this work, like in his inaugural dissertation,\index{doctoral dissertation!Kant}\index{Kant!doctoral dissertation} Kant addresses the fundamental questions that were addressed before him by Leibniz, Newton and others, namely, \emph{What is space? What is time? What is the relation between space, time and the mind? Is this relation real or ideal?  Do space and time have subjective existence, beyond our intuition of them?  Are they empirical concepts? Are they substances or the product of our mind? Do they exist independently of objects and their relation? Are they necessary tools for our understanding?} 
 Kant also analyses our representation of space and its relation to geometry. Elaborating on these most difficult questions is the subject of the fundamental contribution of Kant to philosophy.
  
  Weyl's\index{Weyl, Hermann (1885--1955)} book \emph{Space, time, matter},\footnote{We already recalled that the triad Matter, Space and Time is \emph{par excellence} an Aristotelean theme. There are numerous references regarding this subject, and the best way for the reader to get into this is to skim Aristotle's works. Some of these works are listed in the bibliography, but there are many others.} (first edition 1918),\footnote{The book, under the German title \emph{Raum, Zeit, Materie}, appeared in English translation in 1922.}  is an introduction to the theory of relativity, based on lectures he gave at Zurich's ETH. This work of Weyl is a celebration of the idea that Einstein's theory of relativity is an accomplishment of Riemann's geometry. In the introduction, Weyl writes: 
\begin{quote}\small
It was my wish to present this great subject as an illustration of the intermingling of philosophical, mathematical, and physical thought, a study which is dear to my heart. This could be done only by building up the theory systematically from the foundations, and by restricting attention throughout the principles. But I have not been able to satisfy these self-imposed requirements: the mathematician predominates at the expense of the philosopher.
\end{quote}
The mathematician's role is played essentially by Riemann. 
In Riemannian geometry,\index{Riemannian geometry}\index{geometry!Riemannian} the space (a manifold) is equipped at each tangent space with a quadratic form defining a geometry which is Euclidean.  Weyl comments on this fact and on its relation with physics. He writes (\cite{Weyl-STM} p. 91): 
\begin{quote}\small
The transition from Euclidean geometry to that of Riemann is founded in principle on the same idea as that which led from physics based on action at a distance to physics based on infinitely close action. We find by observation, for example, that the current flowing along a conduction wire is proportional to  the difference of potential between the ends of the wire (Ohm's Law). But we are firmly convinced that this result of measurement applied to a long wire does not represent a physical law in its most general form; we accordingly deduce this law by reducing the measurements obtained to an infinitely small portion of wire. But this means we arrive at the expression on which Maxwell's\index{Maxwell, James Clerk (1831--1879)} theory is founded. Proceeding in the reverse direction, we derive from this differential law by mathematical processes the integral law, which we observe directly, on the supposition that conditions are everywhere similar (homogeneity). We have the same circumstance here. The fundamental fact of Euclidean geometry is that the square of the distance between two points is a quadratic form of the relative co-ordinates of the two points (\emph{Pythagoras Theorem}.) \emph{But if we look upon this law as being strictly valid only for the case when these two points are infinitely near, we enter the domain of Riemann's geometry}. [...] We pass from Euclidean ``finite" geometry to Riemann's ``infinitesimal" geometry in a manner exactly analogous to that by which we pass from ``finite" physics to ``infinitesimal" (or ``contact") physics. 
\end{quote}

Weyl continues (\cite{Weyl-STM} p. 92):

\begin{quote}\small
The principle of gaining knowledge of the external world from the behavior of its infinitesimal parts is the mainspring of the theory of knowledge in infinitesimal physics as in Riemann's geometry, and, indeed, the mainspring of all the eminent work of Riemann, in particular, that dealing with the theory of complex functions. 
\end{quote}

In the same book, Weyl writes (\cite{Weyl-STM} p. 98): 
\begin{quote}\small Riemann rejects the opinion that had prevailed up to his own time, namely, that the metrical structure of space is fixed and inherently independent of the physical phenomena for which it serves as a background,  and that the real content takes possession of it as of residential flats. \emph{He asserts, on the contrary, that space in itself is nothing more than a three-dimensional manifold devoid of all form; it acquires a definite form only through the advent of the material content filling it and determining its metric relations}.
\end{quote}

And then (\cite{Weyl-STM} p. 102):
\begin{quote}\small
Riemann, in the last words of  the above quotation, clearly left the real development of his ideas in the hands of some subsequent scientist whose genius as a physicist could rise to equal flights with his own as a mathematician. After a lapse of seventy years this mission has been fulfilled by Einstein.
\end{quote}

Relativity theory is based on the fact that space and time cannot be separated and form a four-dimensional continuum in one of the senses that Riemann intuited. Einstein  made a profound relation between Riemannian geometry and physics,  in particular in his discovery that gravity is the cause of curvature of physical space. Einstein's equation,\index{Einstein equation}\index{equation!Einstein} published for the first time in 1915, which is the main partial differential equation of general relativity, expresses a relation between energy, gravitation and the curvature of spacetime.\index{Minkowski spacetime}\index{spacetime} In this setting, the Lorentzian metric\index{metric!Lorentzian}\index{Lorentzian metric} encodes the gravitational effects, and the notion of curvature plays a central role.  At several places, Einstein expressed his debt to Riemann.  Let us quote him from \cite{Einstein-Ideas} (p. 281): 
\begin{quote}\small
But the physicists were still far removed from such a way of thinking; space was still, for them, a rigid, homogeneous something, incapable of changing or assuming various states. Only the genius of Riemann, solitary and uncomprehended, had already won its way to a new conception of space, in which space was deprived of its rigidity, and the possibility of its partaking in physical events  was recognized. This intellectual achievement commands our admiration all the more for having preceded Faraday's and Maxwell's\index{Maxwell, James Clerk (1831--1879)} field theory of electricity.
\end{quote}
  
  We end this section by quoting Riemann, and his concerns about physics. In a letter to his father, written February 5, 1852 \cite{Riemann-Letters}, right after he submitted his  Habilitationsschrift,index{Riemann! habilitation text}\index{habilitation text!Riemann} Riemann writes:
\begin{quote}\small
Right after the submission of my Habilitationsschriftindex{Riemann! habilitation text}\index{habilitation text!Riemann} I resumed my investigations into the coherence of the laws of Nature and got so involved in it that I could not tear myself loose. The continuing preoccupation with it has become bad for my health, in fact, right after New Year's my usual affliction set in which such persistence, that I could only obtain relief through the strongest remedies. As a result I felt very ill, felt unable to work, and sought to again put my health in order through long walks.
\end{quote}
On June 26 of the next year, he writes to his brother on the same subject:
\begin{quote}\small
I had completed my habilitation paper at the beginning of December, submitted it to the dean, and soon after once again turned to my investigation on the coherence of the fundamental laws in physics; also that I so immersed myself in it that when the theme for my examination lecture was posted at the colloquium, I could not immediately tear myself away. Rightly after I came down sick, partly, of course, as a result of too much brooding, and partly as a result of sitting a lot in my room during bad weather. 
\end{quote}

\section{The \emph{Commentatio} and the  \emph{Gleichgewicht der Electricit\"at}}

Riemann developed some of his mathematical ideas  introduced in his Habilitationsvortrag in a paper, written in Latin, whose extended title is \emph{Commentatio Mathematica, qua respondere tentatur quaestioni ab Ill${}^{ma}$ Academia Parisiensi propositae: Trouver quel doit \^etre l'\'etat calorifique d'un corps solide homog\`ene ind\'efini pour qu'un syst\`eme de courbes isothermes, \`a un instant donn\'e, restent isothermes apr\`es un temps quelconque, de telle sorte que la temp\'erature d'un point puisse s'exprimer en fonction du temps et de deux autres variables ind\'ependantes} (A mathematical treatise in which an attempt is made to answer the question proposed by the most illustrious Academy of Paris:  To find what must be the thermal state of an indefinite homogeneous solid body so that a system of isothermal curves, at a given instant, remain isothermal after in arbitrary time, in such a way that the temperature at a point can be expressed in terms of time and of two other independent variables). The memoir, as the name indicates, was presented as a contribution to a problem which was proposed for competition by the Paris Academy of Sciences in 1861.
 Part of the \emph{Commentatio} is translated and commented by Spivak in Chapter 4 of Volume II of his \emph{Comprehensive introduction to differential geometry} \cite{Spivak}.

The problem concerns heat conduction, more precisely, the determination of the temperature of a body endowed with a set of given conductivity coefficients.  From the mathematical point of view, it amounts to finding the solution of a partial differential equation -- an evolution equation. The ``solid body" that is referred to in the statement of the problem becomes, in Riemann's context, a Riemannian manifold. At the same time the terms used have a physical significance. 
It is not surprising that Riemann got interested in that problem, which combines geometry and potential theory, two of his favorite subjects. The word ``isothermal" is also reminiscent of the work done by his mentor, Gauss.

While the Habilitationsvortrag\index{Riemann! habilitation lecture}\index{habilitation lecture!Riemann} is practically devoid of any mathematical formulae, the \emph{Commentatio} is full of them. In fact, it is in  the style of the later papers on Riemannian geometry, and in particular those on general relativity, with their debauchery of  indices.\footnote{The expression is due to \'Elie Cartan, from  \cite{Cartan-lecons}, p. VII: 
``The distinguished favor that the absolute differential calculus of Ricci\index{Ricci calculus} and Levi-Civita\index{Levi-Civita calculus} did for us, and will continue to do, should not prevent us of avoiding the over-exclusively formal computations, where the debauchery of indices hides a geometric reality which is often simple. It is this reality which I tried to bring out."}
Riemann's \emph{Commentatio} also contains new tools that are essential to differential geometry. It is in this paper  
that Riemann introduced his 4-entry curvature tensor. The authors of \cite{FR} consider this paper as a ``contribution to the development of what later became known as tensor analysis." As is well known, this topic became an important tool in general relativity.  There is a general agreement now that Riemann's paper contains several results that are usually attributed to Christoffel,\index{Christoffel, Elwin Bruno (1828-1900)} cf. \cite{FR} and \cite{Zund}, and also the idea of Finsler geometry.\index{Finsler geometry}\index{geometry!Finsler}

Let us quote from this paper. Riemann starts his paper by a summary, in which he declares that he will first solve a more general problem: 
\begin{quote}\small
We shall approach the question proposed by the Academy in such a way that we shall first solve a more general question:
the properties of a body which determine the conduction of heat and the distribution of heat within it such that there exists a system of lines which remain isothermal.
\\
From the general solution of this problem we shall distinguish those cases in which the properties vary from those in which the properties remain constant, that is where the body is homogeneous.
\end{quote}
We recall that in the Habilitationsvortrag, the homogeneity property was proved to be equivalent to the curvature being constant.

  The second part of the paper is concerned with the question under equivalence of passing from one quadratic form to another. This is essential in the theory of the transformations that make tensors coordinate-free forms. The reader can find mathematical commentaries on Riemann's memoir in the paper  \cite{FR}.

  The \emph{Commentatio}, like Riemann's Habilitationsvortrag is difficult to read, but this time because of the density of its mathematical content. 
  Riemann's article did not win the prize, probably because some details in the proofs were missing. (In fact, the prize was not awarded.) The authors of \cite{FR} present a certain number of different and conflicting interpretations of the \emph{Commentatio}, a fact which is uncommon for a mathematical paper. This is another indication of how much Riemann's writing are special and cryptic (even today). 
  
The subject of the unfinished paper \emph{Gleichgewicht der Electricit\"at auf Cylindern mit kreisf\"ormigen Querschnitt und parallelen Axen} (On the equilibrium of electricity on cylinders with circular transverse section and whose axes are parallel) (1857) \cite{Riemann-Gleich} by Riemann, published posthumously in the second edition of his Collected works,  is related to the one of the \emph{Commentatio}. It concerns the distribution of electricity or temperature on infinite cylindrical conductors with parallel generatrices. Riemann gives in this paper a solution of the Dirichlet boundary value problem for plane domains. He declares, at the beginning of the paper, that the physical question considered will be solved if the following mathematical question is solved: On a planar connected surface which simply covers the plane and whose boundary may be arbitrary, to determine a function $u$ of the rectangular coordinates $x,y$ satisfying the equation
\[\frac{\partial^2u}{\partial x^2}+\frac{\partial^2u}{\partial y^2}=0\]
and taking arbitrary values on the boundary.  Riemann's solution makes use of Green's theorem and of Abelian integrals. This work is another illustration of the fact that Riemann equates potential theory with the theory of Riemann surfaces.

 \section{Riemann's other papers} 
 We discuss briefly some other papers of Riemann related to our subject. Needless to say, the fact that we pass rapidly through these papers does not mean they are less important than those which we discussed more thoroughly in the previous sections. 
 
Darboux,\index{Darboux, Gaston (1842--1917)} in his famous \emph{Le\c cons sur la th\'eorie g\'en\'erale des surfaces et les applications g\'eom\'etriques du calcul infinit\'esimal} (A course on the general theory of surfaces and the geometric applications of infinitesimal calculus), 1896, \S\,358, discussing the notion of the adjoint equation of a given linear equation, says that the origin of this notion is contained in Riemann's memoir \emph{\"Uber die Fortpflanzung ebener Luftwellen von endlicher Schwingungsweite} (On the propagation of planar air waves that have finite vibration amplitude), \cite{Riemann-Fortpflanzung} 1860; cf.  \cite{Riemann--French} p. 177. He declares that P. du Bois-Reymond, in his work on partial differential equations as well as in a short article he published in T\"ubingen, called the attention of geometers on that memoir by Riemann. He then presents the work. The content is mathematical, with applications to experimental physics. In the introduction, Riemann writes that his research on this subject is in the lineage of the recent work of Helmholtz\index{Helmholtz@von Helmholtz, Hermann Ludwig Ferdinand (1821--1894)} on acoustics. He says that his results, besides their theoretical interest in the theory of the nonlinear partial differential equations which determine the motion of gases, should give the bases for experimental research on the subject.\footnote{\label{foot:R} We note however that in the announcement of this paper, published in the  G\"ottinger Nachrichten, No. 19 (1859), Riemann begins by stating that he does not claim to give any results that are useful in experimental research. At the end of that announcement, he  mentions connections with acoustics but he says that their verification seems to be very hard, the reason being that they either involve very small tone differences, which are not noticeable, or large variations which involve many parameters, therefore the causes cannot be separated. He also talks about applications to meteorology.} He starts in his paper by recalling the physical laws of Boyle, Gay-Lussac and the recent experiments of Regnault, Joule, Thomson\index{Thomson, William (Lord Kelvin) (1824--1907)} and others. About a hundred years later, commenting on the same memoir in the new edition of Riemann's Collected Papers (1990), Peter Lax writes (\cite{Riemann-Gesammelte} p. 807): ``In this paper, Riemann lays the foundations of the theory of propagation of non-linear and linear waves governed by hyperbolic equations. The concepts introduced here -- Riemann invariants, the Riemann initial value problem, jump conditions for nonlinear equations, the Riemann function for linear equations -- are still the basic building blocks of the theory today."  Riemann states in the announcement of the paper (cf. Footnote \ref{foot:R}) that the solution of that problem would help clarifying a perennial debate that involved the mathematicians Challis, Airy,  Stokes, Pretzval, Doppler and Ettinghausen.  Betti\index{Betti, Enrico (1823--1892)} wrote an extensive technical report on that paper, \cite{Betti-1860}.

 We now briefly review some other papers.

 Riemann declares in the introduction to  the paper \emph{Ein Beitrag zu den Untersuchungen \"uber die Bewegung eines fl\"ussigen gleichartigen Ellipsoides} (A contribution to the investigation of the movement of a uniform fluid ellipsoid)  \cite{Riemann-Untersuchungen} that he is continuing the last work of Dirichlet, that this work is surprising and that it opens up a new path for mathematicians which is independent of the original motivation of Dirichlet, which originates in a question on the heavenly bodies.  This paper is also discussed   in a supplement in the new edition of his Collected Works \cite{Riemann-Gesammelte}, \cite{CL}. 
   Riemann's motivation originates in a writing of Newton,\index{Newton, Isaac (1643--1727)} more precisely in his proof of the fact that the spheroidal (rather than spherical)\footnote{Newton, in his \emph{Principia} (1687), expected a flattening of the earth at the poles, of the order of 1/230. The real shape of the earth was another major controversial issue in the seventeenth and eighteenth centuries, and it opposed the English scientists, represented by Newton, to the French, who considered themselves as the heirs of Descartes,\index{Descartes, Ren\'e (1596--1650)} and who were represented by the astronomer Jacques Cassini (1677--1756), the physicist Jean-Jacques d'Ortous de Mairan (1678--1771) and others who pretended on the contrary that the earth was stretched at the poles.  Huygens\index{Huygens, Christiann (1629--1695)} was on the side of Newton.\index{Newton, Isaac (1643--1727)}
Maupertuis\index{Maupertuis, Pierre Louis Moreau de, (1698--1759)} tried to convince the French Academy of Sciences that the theory of Newton concerning the shape of the earth was sound, and he led an expedition to Lapland, whose aim was to measure the length of a meridian. The expedition, which lasted sixteen months, was successful, and it confirmed Newton's\index{Newton, Isaac (1643--1727)} ideas. The mathematician Alexis-Claude Clairaut (1713--1765) and the Swedish astronomer Anders Celsius (1701--1744) were part of the expedition. 
The controversy on the form of the earth gave rise to an extensive literature, in the seventeenth and eighteenth centuries.  In the \emph{Discours pr\'eliminaire} (Preliminary discourse) of the \emph{Encyclop\'edie} (1751), d'Alembert praises Maupertuis\index{Maupertuis, Pierre Louis Moreau de, (1698--1759)} who dared to take side for the English. He writes: ``The first among us who dared to declare openly that he was Newtonian is the author of the \emph{Discours sur la figure des astres} [...]. Maupertuis thought that one could be a good citizen without blindly adopting the physics of one's country; to attack this physics, he needed a courage for which we have to be grateful to him." Voltaire, who contributed in making Newton's\index{Newton, Isaac (1643--1727)} ideas known in France, was among the few major figures on the continent who stood up for the English. He presents these polemics in his famous \emph{Lettres philosophiques} \cite{Voltaire} (1734) (No. XIV): ``A Frenchman who arrives in London, will find philosophy, like everything else, very much changed there. He had left the world a plenum, and he now finds it a vacuum. At Paris, the universe is seen composed of vortices of subtile matter; but nothing like it is seen in London. In France, it is the pressure of the moon that causes the tides; but in England it is the sea that gravitates towards the moon; so that when you think that the moon should make it flood with us, those gentlemen fancy it should be ebb, which very unluckily cannot be proved. [...] At Paris you imagine that the earth is shaped like a melon, or of an oblique figure ; at London it has an oblate one."} form of the earth is due to its rotation. Newton\index{Newton, Isaac (1643--1727)} gave the following formula for a homogeneous body in gravitational equilibrium and small rotation:
 \[m=\frac{5}{4}\epsilon\]
  where  $\epsilon$ is the ellipticity coefficient, equal to the equatorial radius -- polar radius/mean radius,  and $m$ the centrifugal acceleration/mean gravitational acceleration on the surface. The formula was generalized by MacLaurin, who removed the restriction to small rotations. Later works and clarifications are due to Lagrange\index{Lagrange, Joseph-Louis (1736--1813)} and Jacobi. Dirichlet investigated these problems in his 1856/57 lectures on partial differential equations, which were edited in part by Dedekind in 1860.  Chandrasekhar and  Lebowitz, in a commentary on Riemann's paper \cite{Riemann-Untersuchungen} which is published in Riemann's Collected Works edition \cite{Riemann-Gesammelte}, \cite{CL},
 quote Riemann saying: 
 \begin{quote}\small 
 In his posthumous paper, edited for publication by Dedekind, Dirichlet has opened up, in a most remarkable way, an entirely new avenue for investigations on the motion of a self-gravitating homogeneous ellipsoid. The further development of this beautiful discovery has a particular interest to the mathematician even apart from its relevance to the forms of heavenly bodies which initially instigated these investigations.
 \end{quote} We refer the reader to the analysis of Riemann's paper contained in \cite{Riemann-Gesammelte} p. 811-820, where the authors consider this paper to be ``remarkable for the wealth of new results it contains and for the breadth of its comprehension of the entire range of problems. [...] This much neglected paper [...] deserves to be included among the other great papers of Riemann that are well known." In their conclusion, they write: ``A variety of further developments in astronomy and physics have been made possible by the existence of Riemann's work on ellipsoidal figures. [...] The foregoing brief account of developments in the theory of the classical ellipsoids show how Riemann's investigations, after a lapse of some one hundred years, occupy a central place in theoretical astrophysics today."

Let us now say a few words about Riemann's paper \emph{Beitr\"age zur Theorie der durch die Gauss'sche Reihe $F(\alpha,\beta,\gamma,x)$ darstellbaren Functionen}  (Contribution to the theory of functions representable by Gauss's series $F(\alpha,\beta,\gamma,x)$) \cite{Riemann-Beitrage}.  In the introductory part, Riemann announces that in this paper, he investigates the functions representable by Gauss's series using a new method which essentially  applies to any function satisfying a linear differential equation with algebraic coefficients. He also says that the main reasons for his investigations are the many applications of this function in physics and astronomy. In the announcement of that memoir, published in the \emph{G\"ottinger Nachrichten}, No. 1, 1857,  Riemann recalls that Gauss, in studying these functions, was motivated by astronomy.  Riemann's announcement starts with the words: ``This memoir treats a class of functions which are useful to solve various problems in mathematical physics." These functions are still commonly used today in mathematical physics.

Finally, we say a few words on Riemann's paper on minimal surfaces\index{minimal surface} \cite{Riemann-minimal},\footnote{The paper was published posthumously in 1867, and according to Hattendorf,\index{Hattendorf, Karl (1834--1882)} quoted in \cite{Riemann--French} p. 306, it was written around 1860--1861.} \emph{\"Uber die Fl\"ache vom kleinsten Inhalt bei gegebener Begrenzung}  (On surfaces of minimal area, with a given contour). The problem is to find surfaces with minimal area and with fixed boundary. This problem is also related to physics. Again, the mathematical field to which this question belongs can be traced back to the Greeks, namely to works of Archimedes  on isoperimetry and isoepiphany. The specific question of minimal surfaces belongs to the calculus of variations, more precisely the so-called multi-dimensional calculus of variations. In dimension two, one minimizes area or the Dirichlet functional over spaces of surfaces with a given boundary (whereas in the problems of the classical, one-dimensional calculus of variations, one typically minimizes the length, energy, etc. functional on a space of curves joining two given points). It was probably Euler, in 1744, who discovered the first minimal surface, the catenoid, the surface of least area whose boundary consists of two parallel circles in space \cite{Euler-Methodis}. (The name comes from the fact that this is the surface obtained by rotating a catenary around a line.) One year after, Lagrange, who was 19 years old, studied the question of finding the graph of a surface in space with least area with prescribed boundary in the plane. He found a partial differential equation satisfied by such a surface. This was the birth of the so-called Euler--Lagrange equation.\index{Euler--Lagrange equation}\index{equation!Euler-Lagrange} In 1776, Meusnier\footnote{Jean-Baptiste Marie Charles Meusnier de la Place (1754--1793)\index{Meusnier de la Place, Jean-Baptiste Marie Charles (1754--1793)} was a Revolution general, a geometer and an engineer. Together with the mathematicians Gaspard Monge\index{Monge, Gaspard (1746--1818)} and Alexandre-Th\'eophile Vandermonde,\index{Vandermonde, Alexandre-Th\'eophile (1735--1796)} he belonged to the \emph{soci\'et\'e patriotique du Luxembourg}, a patriotic-revolutionary movement. In mathematics, he is known for the Meusnier Theorem on the curvature of surfaces, and for the discovery of the helicoid, a ruled minimal surface.} interpreted Lagrange's\index{Lagrange, Joseph-Louis (1736--1813)} equation as the  vanishing of the mean curvature.  Monge also made substantial contributions to the subject of minimal surfaces. Riemann's contribution (1860--1861) concerns the solution for some given boundary curves. Riemann gave a one-parameter family of examples of minimal surfaces. It was proved recently that the plane, the helicoid, the catenoid and the one-parameter family discovered by Riemann form exactly the set of complete properly embedded, minimal planar domains in $\mathbb{R}^3$, see \cite{MP}.
Weierstrass\index{Weierstrass, Karl (1815--1897)} made the relation between the Euler--Lagrange\index{Euler--Lagrange equation}\index{equation!Euler--Lagrange}  and the Cauchy--Riemann equations. Schwarz obtained results on the same question. The extensive study of minimal surfaces based on soap films was conducted by  Plateau\footnote{Joseph Plateau (1801--1883)\index{Plateau, Joseph (1801--1883)}  obtained a doctorate in mathematics and then became a  physicist. By his experiments on the retina, and for several machines he invented, Plateau is among the first scientists who contributed to the bases of moving images (cinema).} around the year 1873. Besides the relation with soap films, minimal surfaces appear in physics, in particular in hydrodynamics. 
We again cite
Klein, from his article on \emph{Riemann and his significance for the development of modern mathematics} (1895) \cite{Klein1923}: 
\begin{quote}\small
Perhaps less attention has been paid to another physical application in which Riemann's ways of looking at things are laid under contribution in a most attractive manner. I have in mind the theory of \emph{minimum surfaces} [...] the problem is to determine the shape of the least surface that can be spread out in a rigid frame, -- let us say, the form of equilibrium of a fluid lamina that fits in a given contour. It is noteworthy that, with the aid of Riemann's methods, the known functions of analysis are just sufficient to dispose of the more simple cases.  
\end{quote}
This paper on minimal surfaces is analyzed by Yamada in the present volume, \cite{Yamada}.

One could also talk about Riemann's paper on the zeta function, \emph{\"Uber die Anzahl der Primzahlen unter einer gegebenen Gr\"osse} (on the number of prime numbers less than
a given quantity) \cite{Riemann-primes},  recalling that the apparent chaotic distribution of primes has been shown to match the classical random models which describe physical phenomena. 
 \section{Conclusion} 
 
   Beyond Riemann's work which is the subject of the present chapter, one may wonder about the interrelation between mathematics and physics. This subject is complex,  much on it has been said, and adding something new is not a trivial task. Instead, we quote a  text by Picard,\index{Picard, Emile (1856--1941)} from his opening address at the 1920 International Congress of Mathematicians.  Picard is one of the main advocates of the theory of functions of one complex variable, a subject that was dear to Riemann.\footnote{The congress took place at Strasbourg, the place where the idea of the present book came to the editors, at the occasion of two conferences in 2014, the 93th and 94th Encounters between mathematicians and theoretical physicists.  The theme of the 93th encounter was ``Riemann, topology and physics", and that of the 94th was ``Riemann, Einstein and geometry."\index{Einstein, Albert (1879---1955)}.} He writes in \cite{Picard-Strasbourg}: 
  \begin{quote}\small
  Any physical theory which is sufficiently elaborate takes necessarily a mathematical form; it seems that the actions and reactions between spirit and objects gradually brought the formation of moulds where reality could fit, at least in part. For sure, many concepts created by mathematicians did not find yet any application in the study of physical phenomena, but history of science shows that it was reckless to assert that such or such notion would never be used one day. Geometers like to recall the word of the great mathematician Lagrange\index{Lagrange, Joseph-Louis (1736--1813)} who, one day, comparing mathematics to an animal of which every part can be eaten, said: ``Mathematics is like a pig, everything in it is good."\footnote{[Toute th\'eorie physique, suffisamment \'elabor\'ee, prend n\'ecessairement une forme math\'ematique ; il semble que les actions et r\'eactions entre l'esprit et les choses ont amen\'e peu à peu à former des moules où peut, partiellement au moins, s'ins\'erer le r\'eel. Certes, beaucoup de concepts cr\'e\'es par les math\'ematiciens n'ont pas trouv\'e encore d'applications dans l'\'etude des ph\'enom\`enes physiques, mais l'histoire de la science montre qu'il \'etait t\'em\'eraire d'affirmer que telle ou telle notion ne sera pas un jour utilis\'ee. Les g\'eom\`etres aiment à rappeler le mot du grand math\'ematicien Lagrange qui, comparant un jour les math\'ematiques à un animal dont on mange toutes les parties, disait: ``Les math\'ematiques sont comme le porc, tout est bon."]
}
  \end{quote}
%

\noindent \emph{Acknowledgements.--} The author is grateful to Vincent Alberge, Hubert Goenner, Marie Pascale Hautefeuille, Manfred Karbe, Victor Pambuccian and Arkady Plotnitsky  who read  preliminaries versions of this chapter, made comments and suggested corrections.

\printindex

 \end{document}